\documentclass[11pt, oneside]{amsart}   	
\usepackage{amsmath,amsthm,amscd,amsfonts, amssymb, mathrsfs}
\usepackage{geometry}                		
\usepackage{graphicx}				
\usepackage{amssymb}
\newcommand{\sect}[1]{\section{#1}\setcounter{equation}{0}}

\font\mbn=msbm10 scaled \magstep1
\font\mbs=msbm7 scaled \magstep1
\font\mbss=msbm5 scaled \magstep1
\newfam\mbff
\textfont\mbff=\mbn
\scriptfont\mbff=\mbs
\scriptscriptfont\mbff=\mbss   

\newcommand{\Di}      {\mathbb{D}}

\newcommand{\N}       { \mathbb{N}}
\newcommand{\Z}        {\mathbb{Z}  }  
\newcommand\Co           {{\mathbb C}}

\newtheorem{Th}{Theorem}[section]
\newtheorem{Lm}[Th]{Lemma}
\newtheorem{C}[Th]{Corollary}

\newtheorem{Prop}[Th]{Proposition}
\newtheorem{R}[Th]{Remark}

\newtheorem*{Notation}{Notation}

\newtheorem*{Theorem}{Theorem}

\newtheorem*{Th A}{Theorem A}
\newtheorem*{Th B}{Theorem B}
\begin{document}

\title[Runge type Approximation Theorems for $H^\infty$ Maps ]{Dense Stable Rank and Runge Type Approximation Theorems for $\mathbf{H^\infty}$ Maps}
\author{Alexander Brudnyi }
\address{Department of Mathematics and Statistics\newline
\hspace*{1em} University of Calgary\newline
\hspace*{1em} Calgary, Alberta, Canada\newline
\hspace*{1em} T2N 1N4}
\email{albru@math.ucalgary.ca}


\thanks{Research is supported in part by NSERC}

\begin{abstract}
Let $H^\infty(\Di\times\N)$ be the Banach algebra of bounded holomorphic functions defined on the disjoint union of countably many copies of the open unit disk $\Di\subset\Co$.
We show that the dense stable rank of $H^\infty(\Di\times\N)$ is one and using this fact prove some  nonlinear Runge-type approximation theorems for $H^\infty(\Di\times\N)$ maps.  Then we apply these results to obtain a priori uniform estimates of norms of approximating maps in similar approximation problems for algebra $H^\infty(\Di)$.
\end{abstract}

\date{}

\maketitle

\sect{Formulation of Main Results}
\noindent {\bf 1.1.} Let $H^\infty(U)$ denote the Banach algebra of bounded holomorphic functions on a complex manifold $U$ equipped with pointwise multiplication and supremum norm.
In this paper we continue to study properties of the algebra
$H^\infty(\Di\times\N)$, where $\Di\subset\Co$ is the open unit disk, started in \cite{Br2}. 
In particular, we prove that the dense stable rank of $H^\infty(\Di\times\N)$ is one and using this fact prove some nonlinear Runge-type approximation theorems for this algebra. As a consequence, we obtain a priori uniform estimates of norms of approximating maps in similar approximation problems for $H^\infty$ maps defined on subsets of $\Di$.
 
To formulate our results we recall some definitions and notations.

Let $A$ be an associative ring with identity $1$. An element $a=(a_1,\dots, a_n)\in A^n$ is  {\em unimodular} if there exists $b=(b_1,\dots, b_n)\in A^n$ such that $\sum_{i=1}^n b_i a_i=1$. By $U_n(A)\subset A^n$ we denote the set of unimodular elements of $A^n$. An element $a=(a_1,\dots, a_n)\in U_n(A)$ is said to be {\em reducible} if there exist $c_1,\dots, c_{n-1}\in A$ such that
\[
(a_1+c_1 a_n,\dots, a_{n-1}+c_{n-1}a_n)\in U_{n-1}(A).
\]
Ring $A$ is said to have a stable rank at most $n-1$ if every $a\in U_n(A)$ is reducible. The {\em stable rank} of $A$, denoted by ${\rm sr}(A)$, is the least $n-1$ with this property. We write ${\rm sr}(A)=\infty$ if such $n$ does not exist.

The concept of the stable rank introduced by Bass \cite{B}  plays an important role in some stabilization problems of algebraic $K$-theory  analogous to that of dimension in topology.
Despite a simple definition, ${\rm sr}(A)$ is often quite difficult to calculate even for relatively uncomplicated rings $A$ (cf. \cite{V}). An example of such calculation is the  classical Treil theorem asserting that ${\rm sr}(H^\infty)=1$. (We set $H^\infty:=H^\infty(\Di)$.) In fact, Treil proved a much stronger result:
\begin{Theorem}[\cite{T}, Theorem 1]
Let $f,g\in H^\infty$, $\|f\|_{H^\infty}\le 1$, $\|g\|_{H^\infty}\le 1$ and
\begin{equation}\label{treil1}
\inf_{z\in\Di}(|f(z)|+|g(z)|)=:\delta>0.
\end{equation}
Then there exists a function $G\in H^\infty$ such that the function $\Phi=f+gG$ is invertible in $H^\infty$, and moreover $\|G\|_{H^\infty}\le C$, $\|\Phi^{-1}\|_{H^\infty}\le C$, where the constant $C$ depends only on $\delta$.
\end{Theorem}
\noindent (Here and below for a normed space $B$ its norm is denoted by $\|\cdot\|_B$.)

It was observed in  \cite[Thm.\,2]{To} that uniform estimates of norms of $G$ and $\Phi$ of the theorem imply that ${\rm sr}(H^\infty(\mathbb D\times\N))=1$ and de facto the latter statement is equivalent to the Treil theorem. In the present paper, we show that some other results for  $H^\infty(\mathbb D\times\N)$ lead to analogs of the above theorem with estimates of norms of the corresponding `solutions' depending only on numerical data. Such estimates are of importance in many nonlinear problems for algebra $H^\infty$, see, e.g., \cite{Ga} for the references.\smallskip

 It is well known \cite{CL}, \cite{CS1} that if $A$ is a commutative complex unital Banach algebra, then 
${\rm sr}(A)\le n$ iff for every epimorphism of $A$ onto a commutative complex unital Banach algebra $B$ the induced map $U_n(A)\rightarrow U_n(B)$ is onto. This leads to a more general notion of stability, see \cite{CS2}:

The {\em dense stable rank} of a commutative complex unital Banach algebra $A$, denoted by ${\rm dsr}(A)$, is the least integer $n$ such that for every dense image morphism of $A$ in a commutative complex unital Banach algebra $B$  the induced map $U_n(A)\rightarrow U_n(B)$ has dense image. (As before, ${\rm dsr}(A)=\infty$ if there is no such $n$.)   Clearly, ${\rm sr}(A)\le {\rm dsr}(A)$; the question whether equality holds is open.

 Using some arguments of the proof of \cite[Thm.\,1]{T}, Su\'arez \cite[Cor.\,4.6]{S1} proved that ${\rm dsr}(H^\infty(\Di))=1$. Based on this result we prove:
 \begin{Th}\label{te1}
${\rm dsr}(H^\infty(\mathbb D\times\N))=1$.
\end{Th}
\begin{Notation}
{\rm 
In what follows, $\mathfrak M(A)$  stands for the maximal ideal space of a commutative complex unital Banach algebra $A$ (i.e., the set of nonzero homomorphisms $A\rightarrow\mathbb C$ equipped with the {\em Gelfand topology}). Also, if $X$ and $Y$ are topological spaces, then $[X]$ denotes the set of connectivity components of $X$ and $[X,Y]$ denotes the set of homotopy classes of continuous maps from $X$ to $Y$.}
\end{Notation}
Theorem \ref{te1} asserts that each dense image morphism of complex unital Banach algebras $\varphi: H^\infty(\mathbb D\times\N)\rightarrow A$  induces a dense image map  $U_n(H^\infty(\mathbb D\times\N))\rightarrow U_n(A)$, $n\in\N$, see \cite[Prop.\,3.2]{CS2}. Hence, it induces a surjective map
$[U_n(H^\infty(\mathbb D\times\N))]\rightarrow [U_n(A)]$. For $n=1$ due to the Arens-Royden theorem the latter implies that $\varphi$ induces an epimorphism of \v{C}ech cohomology groups 
$H^1(\mathfrak M(H^\infty(\mathbb D\times\N)),\Z)\rightarrow H^1(\mathfrak M(A),\Z)$.  In fact, this statement is equivalent to Theorem \ref{te1} as well as to the next `quantitative' Runge-type approximation theorem similar to \cite[Thm.\,1]{T}.  In the  formulation of the theorem we use the following notation.

In the sequel, a set of the form
\begin{equation}\label{e3.1}
\Pi_{c;\nu}^n:=\left\{z\in\Di\, :\, \max_{1\le i\le n}|f_i(z)|<\nu,\, f_i\in H^\infty,\, \|f_i\|_{H^\infty}\le c,\, 1\le i\le n\right\}
\end{equation}
is called an open $n$-polyhedron. We set $\Pi_{c}^n:=\Pi_{c;1}^n$. 

 \noindent The group of invertible elements of a unital algebra $A$ is denoted by $A^{-1}$.

\noindent We write $C=C(\alpha_1,\alpha_2,\dots)$ if the constant $C$ depends only on $\alpha_1,\alpha_2,\dots$.

\begin{Th}\label{lem3.1}
Suppose $g\in (H^\infty(\Pi_c^n))^{-1}$ satisfies 
\begin{equation}\label{e3.2}
\|g^{\pm 1}\|_{H^\infty(\Pi_c^n)}\le M.
\end{equation}
Given $0<\delta<1$ and $\varepsilon>0$, there exist a constant $C=C(M,c,  n,\delta,\varepsilon)$ and a function $h\in (H^\infty)^{-1}$ such that
 \begin{equation}\label{e3.3a}
\|h^{\pm 1}\|_{H^\infty}\le C\quad {\rm and}\quad
\bigl\|(h-g)|_{\Pi_{c;\delta}^n}\bigr\|_{H^\infty(\Pi_{c;\delta}^n)}\le\varepsilon. 
\end{equation}
\end{Th}

\medskip 

\noindent {\bf 1.2.} In this part, we present some extensions of Theorems \ref{te1} and \ref{lem3.1}. To this end, we recall some concepts of the Novodvorski--Taylor theory \cite{Ta}.\smallskip

Let $A$ be a commutative complex unital Banach algebra. By definition, $\mathfrak M(A)$ is a  weak$^*$ compact subset of the dual space $A^*$.
If $U\subset A^*$ is a weak$^*$ open neighbourhood of $\mathfrak M(A)$, we denote by  $\mathcal O(U,\Co^n)$ the vector space of weak$^*$ continuous holomorphic maps of $U$ in $\Co^n$. (Recall that a continuous in the norm topology on $A^*$ map $U\rightarrow\Co^n$ is {\em holomorphic} if all its complex directional derivatives exist.) We equip 
$\mathcal O(U,\Co^n)$ with the topology of uniform convergence on weak$^*$ compact subsets of $U$. Then we define $\mathcal O(\mathfrak M(A),\Co^n)$  to be the inductive limit of spaces $\mathcal O(U,\Co^n)$ as $U$ ranges over all  weak$^*$ open neighbourhoods of $\mathfrak M(A)$ equipped with the inductive limit topology. We regard $\mathcal O(\mathfrak M(A),\Co^n)$ as a topological algebra with the product induced by
coordinatewise multiplication of maps in $\mathcal O(U,\Co^n)$. Similarly, we consider $A^n$ as the Banach algebra equipped with coordinatewise multiplication. Then there exists a continuous algebra epimorphism, the {\em holomorphic functional calculus}, $T_A^n: \mathcal O(\mathfrak M(A),\Co^n)\rightarrow A^n$
such that 
\begin{itemize}
\item[(i)]
The composition $(\,\hat{\, }^{\,n})\circ T_A^n: \mathcal O(\mathfrak M(A),\Co^n)\rightarrow C(\mathfrak M(A),\Co^n)$, where 
$\, \hat{\, }^{\, n} : A^n\rightarrow C(\mathfrak M(A),\Co^n)$, 
\[
(\, \hat{\, }^{\, n}) (a_1,\dots,a_n)(\xi):=(\xi(a_1),\dots,\xi(a_n)),\quad \xi\in\mathfrak M(A),
\] 
is the $n$-fold product of the {\em Gelfand transform}, is the restriction map $f\mapsto f|_{\mathfrak M(A)}$;\medskip
\item[(ii)] 
The composition of $T_A^n$ and the $n$-fold product of the injective map $A\rightarrow \mathcal O(\mathfrak M(A),\Co)$ defined by the natural embedding $A\hookrightarrow A^{**}$
is ${\rm id}_{A^n}$, the identity map of $A^n$.
\end{itemize}

Next, for a complex submanifold $\mathcal M\subset\Co^n$ we define $\mathcal O(U,\mathcal M)\subset\mathcal O(U,\Co^n)$ as the subset of maps with images in $\mathcal M$. Applying to the family of spaces $\mathcal O(U,\mathcal M)$ the inductive limit construction we get a subspace $\mathcal O(\mathfrak M(A),\mathcal M)\subset\mathcal O(\mathfrak M(A),\Co^n)$. Then we set
\begin{equation}\label{eq1.4}
A_{\mathcal M}:= T_A^n(\mathcal O(\mathfrak M(A),\mathcal M))\, (\subset A^n).
\end{equation}

It is known that if $\mathcal M\subset\Co^n$ is a complex submanifold and a homogeneous manifold (i.e., it is equipped with a holomorphic and transitive action of a complex Lie group, see, e.g., \cite{A} for the references), then $A_{\mathcal M}$ is
locally path-connected and the map $[A_\mathcal M]\rightarrow [\mathfrak M(A),\mathcal M]$ induced by $\, \hat{\, }^{\, n}$ is a bijection, see \cite[Sec.\,2.7]{Ta}. In fact, the same is true for the class of manifolds $\mathcal M$ subject to the following definition.\smallskip

\noindent $(\triangledown)$ A complex manifold $\mathcal M$ is said to be {\em Oka} if every holomorphic map $f: K\rightarrow \mathcal M$ from a neighbourhood of  a compact convex set $K \subset\Co^k$, $k\in\N$, can be approximated uniformly on $K$ by entire maps $\Co^k\rightarrow \mathcal M$.

The class of Oka manifolds includes, in particular, complex homogeneous manifolds, complements in $\Co^k$, $k>1$, of complex algebraic subvarieties of codimension $\ge 2$ and of compact polynomially convex sets, Hopf manifolds (i.e., nonramified holomorphic quotients of $\Co^n\setminus\{0\}$). Also,  holomorphic fibre bundles whose bases and fibres are Oka manifolds are Oka manifolds as well. (We refer to the book \cite{F} and the paper \cite{K} for other examples and basic results of the theory of Oka manifolds.)

\begin{Prop}\label{prop1.3}
 Let $\mathcal M\subset\Co^n$ be a complex submanifold and an Oka manifold. Then $A_{\mathcal M}$ is
locally path-connected and the map $[A_\mathcal M]\rightarrow [\mathfrak M(A),\mathcal M]$ induced by $\, \hat{\, }^{\, n}$ is a bijection.
\end{Prop}

Now we move on to the formulation of the nonlinear Runge-type approximation theorem for commutative complex unital Banach algebras.

Let $r: A_1\rightarrow A_2$ be a dense image morphism of commutative complex unital Banach algebras. Then the transpose map $r^*$ embeds $\mathfrak M(A_2)$ into $\mathfrak M(A_1)$. In turn, the pullback with respect to $r^*|_{\mathfrak M(A_2)}$ denoted by $r_c$ maps $C(\mathfrak M(A_1))$ surjectively onto $C(\mathfrak M(A_2))$.

Let $\mathcal M\subset\Co^n$ be a complex submanifold. Then the $n$-fold product of $r$ (denoted by $r^n$) maps $(A_{1})_\mathcal M$ in $(A_2)_\mathcal M$ (see Proposition \ref{prop4.0} below). Let $\overline{r^n((A_1)_\mathcal M)}$ be the closure of $r^n((A_1)_\mathcal M)$ in $(A_2)_\mathcal M$.
\begin{Th}\label{te4.1}
Assume that $\mathcal M$ is an Oka manifold. The following is true:
\begin{itemize}
\item[(1)]
\[
\overline{r^n((A_1)_\mathcal M)}
=\bigl\{a\in (A_2)_{\mathcal M}\, :\ \hat{\, }^{\, n}(a)\in {\rm range}\, (r_c)^n\bigr\}.
\]
\item[(2)] If $r$ is an epimorphism, then $r^n((A_1)_\mathcal M)$ is a closed subset of $(A_2)_\mathcal M$.
\end{itemize}
\end{Th}
Thus, identifying $\mathfrak M(A_2)$ with its image under $r^*$ in $\mathfrak M(A_1)$ we obtain that
an element $a\in (A_2)_{\mathcal M}$ is approximated by images under $r^n$ of elements from $(A_1)_{\mathcal M}$ if and only if its Gelfand transform $\,\hat{\, }^{\, n}(a)\in C(\mathfrak M(A_2),\mathcal M)$ is extended to a map from $C(\mathfrak M(A_1),\mathcal M)$. 
\smallskip
\begin{R}\label{rem1.5}
{\rm
Let $A$ be a commutative complex unital Banach algebra.
We set for a complex submanifold $\mathcal M\subset\Co^n$,
\begin{equation}\label{eq1.3}
A^\mathcal M:=  \{ (a_1,\dots, a_n)\in A^n\, :\, (\xi(a_1),\dots, \xi(a_n))\in \mathcal M\quad \forall\, \xi\in\mathfrak M(A)\}.
\end{equation}
Clearly, $A_\mathcal M \subseteq A^\mathcal M$ and these sets coincide if either $\mathcal M$ is an open subset of $\Co^n$ or $A$ is semisimple (i.e., the Gelfand transform $\,\hat{\,}:A\rightarrow C(\mathfrak M(A))$ is injective), see \cite[Sec.\,2.8]{Ta}.

Let $R_A\subset A$ be the Jacobson radical, i.e., the kernel of the Gelfand transform of $A$. Then $j_A:A\rightarrow A/R_A=:A_s$ is an epimorphism of  complex unital Banach algebras and $A_s$ is semisimple. Moreover, the transpose map $(j_A)^{*}$ maps $\mathfrak M(A_s)$ homeomorphically onto $\mathfrak M(A)$. Applying Theorem \ref{te4.1} to $A_1:=A$, $A_2:=A_s$ we obtain:\smallskip

\noindent (A)} Suppose $\mathcal M\subset\Co^n$ is a complex submanifold and an Oka manifold. Then  $(j_A)^n$ maps $A_\mathcal M$ surjectively onto $(A_s)^{\mathcal M}$.
Moreover, the map $[A_{\mathcal M}]\rightarrow \left[(A_s)^{\mathcal M}\right]$ of sets of connectivity components induced by  $(j_A)^n$ is a bijection.
\smallskip

{\rm Let $r: A_1\rightarrow A_2$ be a dense image morphism of commutative complex unital Banach algebras. Clearly, $r(R_{A_1})\subset R_{A_2}$. Hence, there is a dense image morphism $r_s: (A_{1})_s\rightarrow (A_2)_s$ such that $j_{A_2}\circ r=j_{A_1}\circ r_s$. Now Theorem \ref{te4.1} leads to the following statement.
\smallskip

\noindent (B)} Suppose that $\mathcal M\subset\Co^n$ is a complex submanifold and an Oka manifold. An element $a\in (A_2)_\mathcal M$ is approximated by  images under $r^n$ of elements from $(A_1)_{\mathcal M}$ if and only if $j_{A_2}(a)\in ((A_2)_s)^\mathcal M$ is  approximated by images under $(r_s)^n$ of elements from $((A_1)_s)^{\mathcal M}$.\smallskip

{\rm This reduces the general Runge-type theorem for commutative complex unital Banach algebras to the case of semisimple algebras.
}
\end{R}

As a corollary of Theorem \ref{te4.1} we obtain an extension of Theorem \ref{te1}.
To formulate the result, we recall the following definition.

\smallskip

\noindent $(\triangledown)$ A path-connected topological space $X$ is  $i$-{\em simple} if for each $x\in X$ the fundamental group $\pi_1(X,x_0)$ acts trivially on the $i$-homotopy group $\pi_i(X,x)$ (see, e.g., \cite[Ch.\,IV.16]{Hu3} for the corresponding definitions and results).  

For instance, $X$ is $i$-simple if group $\pi_i(X)=0$ and
$1$-simple if and only if group $\pi_1(X)$ is abelian. Also, every path-connected topological group is $i$-simple for all $i$. The same is true for a complex manifold biholomorphic to the quotient of a connected complex Lie group by a {\em connected} closed Lie subgroup, see, e.g., \cite[(3.2)]{Hu1}.\smallskip

Let $\varphi: H^\infty(\mathbb D\times\N)\rightarrow A$ be a dense image morphism of complex unital Banach algebras and  let $\mathcal M\subset\Co^n$ be a connected submanifold and an Oka manifold. 
\begin{Th}\label{teo1.3}
Assume that a finite unbranched covering of $\mathcal M$ is $i$-simple for $i=1,2$.
Then the image of  $(H^\infty(\mathbb D\times\N))^\mathcal M$ under $\varphi^n$ is a dense subset of $A_\mathcal M$. 

\noindent If, in addition, $\varphi$ is an epimorphism, then this image coincides with $A_\mathcal M$.
\end{Th}

Since each connected component of a topological space is open, Theorem \ref{teo1.3}
implies the following result.

\begin{C}\label{teo1.2}
Under hypotheses of Theorem \ref{teo1.3}, the map $\bigl[(H^\infty(\mathbb D\times\N))^\mathcal M\bigr]\rightarrow [A_\mathcal M]$ induced by $\varphi^n$ is a surjection.
\end{C}
\begin{R}
{\rm (1) As was mentioned in Section~1.1, Theorem \ref{te1} follows from Corollary \ref{teo1.2} with $n=1$ and $\mathcal M=\Co^*:=\Co\setminus\{0\}$. \smallskip

\noindent (2) Let $\mathscr O$ be the class of connected Oka manifolds $\mathcal M$ embeddable as complex submanifolds into complex Euclidean spaces and having $i$-simple for $i=1,2$ finite unbranched coverings. This class contains e.g., complements in $\Co^k$, $k>1$, of complex algebraic subvarieties of codimension $\ge 2$ and of compact polynomially convex sets (these manifolds are simply connected, see \cite{F1}), connected Stein Lie groups, quotients of connected reductive complex Lie groups by Zariski closed subgroups (these manifolds are quasi-affine algebraic, see, e.g., \cite[Thm.\,5.6]{A} for the references;  they have $i$-simple finite unbranched coverings because Zariski closed subgroups have finitely many connected components and quotients of connected complex Lie groups by connected closed Lie subgroups are $i$-simple for all $i$, see, e.g., \cite[(3.2)]{Hu1}). Also, direct products of manifolds from class $\mathscr O$ belongs to $\mathscr O$, etc.
}
\end{R}

In the next two results, $\mathcal M\subset\Co^n$ is a connected complex submanifold from class
$ \mathscr O$, $K\subset\mathcal M$ is a compact subset, and $\Pi_c^k\subset\Di$ is an open $k$-polyhedron, see \eqref{e3.1}.

As a corollary of Theorem \ref{teo1.3} we obtain the following extension of Theorem \ref{lem3.1}.
\begin{Th}\label{teo1.5}
Suppose $g\in \mathcal O(\Pi_c^k,\mathcal M)$ is 
 such that $g(\Pi_c^k)\subset K$. Given $0<\delta<1$ and $\varepsilon>0$ 
 there exist a constant $C=C(\mathcal M, K,n, c, k,\delta,\varepsilon)>0$  and a map $h\in (H^\infty)^\mathcal M$
 such that 
 \begin{equation}\label{equ1.6}
\|h\|_{(H^\infty)^n}\le C\quad {\rm and}\quad \|(h-g)|_{\Pi_{c;\delta}^k}\|_{(H^\infty(\Pi_{c;\delta}^k))^n}\le\varepsilon .
\end{equation}
\end{Th}

Further, if $\Pi_{c}^k$ is determined by functions $f_1,\dots, f_k\in H^\infty(\Di)$, then we define an open set $\hat\Pi_{c;\nu}^k\subset\mathfrak M(H^\infty)$ by the formula
\begin{equation}\label{polyh1}
\hat\Pi_{c;\nu}^k:=\left\{\xi\in \mathfrak M(H^\infty)\, :\, \max_{1\le i\le k}|\hat f_i(\xi)|<\nu\right\},\quad \hat\Pi_{c}^k:=\hat\Pi_{c;1}^k.
\end{equation}
Clearly, $\Pi_{c}^k=\hat\Pi_{c}^k\cap\Di$ and by the corona theorem
$\Pi_c^k$ is an open dense subset of $\hat\Pi_{c}^k$. Moreover, due to \cite[Thm.\,3.2]{S1} each $f\in H^\infty(\Pi_c^k)$ admits an extension $\hat f\in C(\hat\Pi_c^k)$.\smallskip

Let $J\subsetneq H^\infty(\Di)$ be a closed ideal and let
\begin{equation}\label{hull1}
{\rm hull}\, J:=\{\xi\in \mathfrak M(H^\infty)\, :\, \hat{f}(\xi)=0\quad \forall\, f\in J\}.
\end{equation}

Let $\Pi_c^k$  be a polyhedron such that 
\begin{equation}\label{subset}
 {\rm hull}\, J \subset \hat\Pi_{c;\delta}^k\quad {\rm for\ some}\quad 0<\delta<1.
\end{equation}
(Such polyhedra always exist, see Remark \ref{rem1.11} below.)
\begin{Th}\label{teo1.9}
Let $g\in (H^\infty)^n$ be such that
\begin{equation}\label{equ1.11}
\|g\|_{(H^\infty)^n}\le b\quad {\rm and}\quad (\,\hat{\,}^{\, n}(g))(\xi)\in K\quad \forall\, \xi\in {\rm hull}\, J.
\end{equation}
Suppose that there is $f\in\mathcal O(\Pi_c^k,\mathcal M)$, $f(\Pi_c^k)\subset K$, such that
\[
\hat f|_{{\rm hull}\, J}=\hat g|_{{\rm hull}\, J}.
\]
Then there exist a constant $C=C(\mathcal M, K,n,b, c, k,\delta)>0$ and a map $h\in (H^\infty)^\mathcal M$ such that
\begin{equation}\label{eq1.12}
\|h\|_{(H^\infty)^n}\le C\quad {\rm and}\quad \,\hat{\,}^{\, n}(h)|_{{\rm hull}\, J}=\,\hat{\,}^{\, n}(g)|_{{\rm hull}\, J}.
\end{equation}
\end{Th}
 
A particular case of Theorem \ref{teo1.9} for $\mathcal M=\Co^*$  follows from Treil's theorem \cite{T} via its cohomological interpretation due to Su\'{a}rez \cite[Thm.\,1.3]{S1} and the Arens-Royden theorem.
\begin{R}\label{rem1.11}
{\rm A compact subset $S\subset\mathfrak M(H^\infty)$ is called {\em holomorphically convex} if for each $\xi\not\in S$ there is $f\in H^\infty$ such that
\begin{equation}\label{eq1.13}
\max_{S}|\hat f|< |\hat f(\xi)|.
\end{equation}
If $U$ is an open neighbourhood of such $S$, then there exists an open polyhedron $\Pi_{c}^k\subset\Di$ such that $S\subset \hat{\Pi}_c^k\subset U$, see, e.g., \cite[Lm.\,5.1]{Br1}. Hence, Theorems \ref{teo1.5} and \ref{teo1.9} imply the following result.

Let $S\subset\mathfrak M(H^\infty)$ be holomorphically convex, $U\subset \mathfrak M(H^\infty)$ be an open neighbourhood of $S$ and  $K\subset\mathcal M$ be as in the above theorems.}
\begin{C}\label{cor1.11}
{\rm (1)} Suppose $g\in C(U,K)$ is such that $g|_{U\cap\Di}$ is holomorphic.
Given $0<\delta<1$ and $\varepsilon>0$ 
 there exist a constant $c=c(\mathcal M, K,n,S,U,\delta,\varepsilon)>0$  and a map $h\in (H^\infty)^\mathcal M$
 such that 
 \begin{equation}\label{eq1.14}
\|h\|_{(H^\infty)^n}\le c\quad {\rm and}\quad \|(\,\hat{\,}^{\, n}(h)-g)|_S\|_{C(S,K)}\le\varepsilon .
\end{equation}
{\rm (2)} Let $S={\rm hull}\, J$ for some closed ideal $J\subsetneq H^\infty$ and let $g\in (H^\infty)^n$ satisfy
\begin{equation}\label{equ1.15}
\|g\|_{(H^\infty)^n}\le b\quad {\rm and}\quad (\,\hat{\,}^{\, n}(g))(\xi)\in K\quad \forall\, \xi\in S.
\end{equation}
Suppose that there is $f\in C(U,K)$ such that $f|_{U\cap\Di}$ is holomorphic and
\[
f|_{S}=\hat g|_{S}.
\]
Then there exist a constant $c=c(\mathcal M, K,n,b, S,U)>0$ and a map $h\in (H^\infty)^\mathcal M$ such that
\begin{equation}\label{eq1.16}
\|h\|_{(H^\infty)^n}\le c\quad {\rm and}\quad \,\hat{\,}^{\, n}(h)|_{{\rm hull}\, J}=\,\hat{\,}^{\, n}(g)|_{{\rm hull}\, J}.
\end{equation}
\end{C}
{\rm The previous result partially extends Theorems 3.7 and 3.8 of  \cite{Br3} and give a priori uniform estimates of norms of approximating and interpolating maps there.}
\end{R}

 \sect{Auxiliary Results}
\noindent  {\bf 2.1.}  Due to the Carleson corona theorem $\Di$ can be regarded as a dense open subset of the maximal ideal space $\mathfrak M(H^\infty)$ equipped with the Gelfand topology. Then the Gelfand transform extends each  $f\in H^\infty$ to a unique continuous function $\hat{f}$ on $\mathfrak M(H^\infty)$. For an open subset $U\subset \mathfrak M(H^\infty)$ the corona theorem implies that $U\cap\Di$ is an open dense subset of $U$. A function in $C(U)$ is called holomorphic if its restriction to $U\cap\Di$ is holomorphic in the standard sense. The set of such functions is denoted by $\mathcal O(U)$. Each $f\in H^\infty(U\cap\Di)$ extends (uniquely) to a bounded function $\hat{f}\in \mathcal O(U)$, see  \cite[Thm.\,3.2]{S1}.  
A compact subset $K\subset U$ is called $\mathcal O(U)$-convex if for each $x\not\in K$ there is $f\in\mathcal O(U)$ such that
\begin{equation}\label{eq2.5}
\max_K |f|<|f(x)|.
\end{equation}
{\bf 2.2.} In what follows, $\Di_r(c):=\{z\in\Co\, :\, |z-c|<r\}$, $r>0$, $c\in\Co$, i.e., $\Di:=\Di_1(0)$. For a subset $S$ of a topological space we denote by  $\bar S$  its closure.
  
Let 
\begin{equation}\label{equ3.1}
\Omega:=\bigl(\Di\cup\Di_1({\scriptstyle \frac 3 2})\bigr)\setminus\bigl\{\mbox{$ \frac 3 2$}\bigr\}.
\end{equation}
The fundamental group $\pi_1(\Omega)$ of $\Omega$ is isomorphic to $\Z$, i.e., $\pi_1(\Omega)=\{a^n\}_{n\in\Z}$ for some $a\in \pi_1(\Omega)$.
Let $r:\Di\rightarrow\Omega$ be the universal covering of $\Omega$. The deck transformation group  $\pi_1(\Omega)$ acts discretely on $\mathbb D$ by M\"{o}bius transformations. Since $r\in H^\infty$, it extends to a function $\hat r\in \mathcal O(\mathfrak M(H^\infty))$ such that $\hat r(\mathfrak M(H^\infty))=\bar\Omega$.  Let $U:=r^{-1}(\Di)\subset\Di$. Since each loop in $\Di$ is contractible in $\Omega$, 
\begin{equation}\label{equ3.1a}
U=\bigsqcup_{g\in\pi_1(\Omega)}g(U')
\end{equation}
for some $U'\subset\Di$ biholomorphic to $\Di$ via $r$. In particular, the map  $s: U\rightarrow \Di\times\Z$, 
\begin{equation}\label{equ3.2}
s(z):=(r(z),n),\quad  z\in a^n(U'),\quad  n\in\Z,
\end{equation}
is biholomorphic.

By the corona theorem $U$ is dense in 
\begin{equation}\label{equ3.4}
\widetilde U:=\left\{x\in \mathfrak M(H^\infty)\, :\, |\hat r(x)|<1\right\}.
\end{equation}
and due to  \cite[Thm.\,3.2]{S1} each $f\in H^\infty(U)$ extends to a (unique) $\hat f\in \mathcal O(\widetilde U)$.

Let us consider $h:={\rm Re}\, r$
and its extension $\hat h:={\rm Re}\, \hat r\in C(\mathfrak M(H^\infty))$. 
Clearly,
\begin{equation}\label{eq2.6}
N:=\left\{x\in\mathfrak M(H^\infty)\, :\, \hat h(x)\le 0\right\}
\end{equation}
is a $\mathcal O(\mathfrak M(H^\infty))$-convex compact subset of  $\widetilde U$.
\begin{Lm}\label{lemma2.2}
Every $\mathcal O(\widetilde U)$-convex subset $K\subset N$  is $\mathcal O(\mathfrak M(H^\infty))$-convex.
\end{Lm}
\begin{proof}
Let $x\not\in K$ and $f\in\mathcal O(\widetilde U)$ satisfy \eqref{eq2.5}. Clearly, we can assume that $x\in N$.
We set
\[
W:=\left\{x\in\mathfrak M(H^\infty)\, :\, \hat h(x)< \frac 12\right\}.
\]
Then $W\Subset\widetilde U$ is an open neighbourhood of $N$. Since $N$ is $\mathcal O(\mathfrak M(H^\infty))$-convex, $f|_W$ can be uniformly approximated on $N$ by functions in $\mathcal O(\mathfrak M(H^\infty))$, see \cite[Cor.\,2.6]{S2} or \cite[Thm.\,1.7]{Br1}. Hence, there is $f'\in \mathcal O(\mathfrak M(H^\infty))$ sufficiently close to $f|_N$ such that \eqref{eq2.5} holds for $f$ replaced by $f'$.
This proves the lemma. 
\end{proof}
\noindent {\bf 2.3.} Let $A(\Di)$ be the disk algebra of holomorphic functions in $\Di$ continuous on $\bar{\Di}$. A compact set $K\subset\Co$ is called polynomially convex if for each $z\not\in K$ there is a holomorphic polynomial $p$ on $\Co$ such that
\[
\max_K |p|<|p(z)|.
\]
\begin{Lm}\label{lem2.2}
Let $K\subset\bar{\Di}$ be a proper polynomially convex set. Then there is a M\"{o}bius transformation $g:\Di\rightarrow\Di$ such that 
\[
g(K)\subset D_{-}:=\{z\in\bar{\Di}\, :\, {\rm Re}\,z\le 0\}.
\]
\end{Lm}
\begin{proof}
Since $K$ is proper and polynomially convex, there exist a point $c\in\mathbb S:=\bar{\Di}\setminus\Di$ and an open disk $\Di_r(c)$ such that $\Di_r(c)\cap K=\emptyset $. We choose some $t\in (1-r,1)$ such that the disk $\Di_{r'}(c')$ of radius $r':=\frac{1}{2}\bigl(\frac{1}{t\bar{c}}-tc\bigr)$ centered at $c':=\frac{1}{2}\bigl(\frac{1}{t\bar{c}}+tc\bigr)$ satisfies
\[
\Di_{r'}(c')\cap\bar{\Di}\subset\Di_r(c).
\]
By the definition, the part of the boundary of $\Di_{r'}(c')$ in $\Di$ is a geodesic in the Poincar\'{e} metric on $\Di$ passing through $tc$. Then the M\"{o}bius transformation
\[
g(z):=\bar{c}\cdot\frac{z-tc}{1-t\bar{c}z},\quad z\in\bar{\Di},
\]
maps this geodesic to the interval $\{{\rm Re}\,z=0\}\cap\Di$ and $\bar{\Di}\setminus \Di_{r'}(c')$ to $D_{-}$. Hence, $g(K)\subset D_{-}$ as required.
\end{proof}
\sect{Proofs of Theorems \ref{te1} and \ref{lem3.1}}
\begin{proof}[{\bf 3.1.} Proof of Theorem \ref{lem3.1}] 
Recall that
\begin{equation}\label{e3.1a}
\Pi_{c;\delta}^n:=\left\{z\in\Di\, :\, \max_{1\le i\le n}|f_i(z)|<\delta,\, f_i\in H^\infty,\, \|f_i\|_{H^\infty}\le c,\, 1\le i\le n\right\},\ \,  \Pi_{c}^n:=\Pi_{c;1}^n.
\end{equation}
We also use notation $\Pi_{c;\delta}^n[F]$ and $\Pi_c^n [F]$, where $F:=\{f_i\}_{1\le i\le n}$, to emphasize dependence of the polyhedron on the family of functions determining it.

We have to prove that given $g\in (H^\infty(\Pi_c^n))^{-1}$ satisfying
\begin{equation}\label{e3.0}
\|g^{\pm 1}\|_{H^\infty(\Pi_c^n)}\le M
\end{equation}
and every $\varepsilon>0$ there exist a constant $C=C(c, M, n,\delta,\varepsilon)$ and a function $h\in (H^\infty)^{-1}$ such that
 \begin{equation}\label{e3.3}
\|h^{\pm 1}\|_{H^\infty}\le C\quad {\rm and}\quad
\bigl\|(h-g)|_{\Pi_{c;\delta}^n}\bigr\|_{H^\infty(\Pi_{c;\delta}^n)}\le\varepsilon. 
\end{equation}

It suffices to prove the lemma for $n$-polyhedra $\Pi_c^n$ determined by functions from the disk algebra $A(\Di)$. Then a standard normal family argument produces the required result in the general case. Clearly, we may assume that closures of such polyhedra are proper subsets of $\bar{\Di}$ (for otherwise, the statement is trivial). In turn, applying Lemma \ref{lem2.2} we may assume that $\Pi_c^n\subset D_{-}$.  If $\Pi_c^n$ is determined by  $F=\{f_i\}_{1\le i\le n}\subset H^\infty$,  then the set $\mathcal D=\{F, c, M, n,\delta,\varepsilon\}$ is called the {\em data}.

According to \cite[Cor.\,2.6]{S2} each function in $H^\infty(\Pi_c^n)$ can be uniformly approximated on subsets $\Pi_{c;\nu}^n$, $0<\nu<1$, by functions from $H^\infty$. Then according  to
\cite[Cor.\,4.6]{S1} applied to the dense morphism $H^\infty\rightarrow A$, $f\mapsto f|_{\Pi_{c;\delta}^n}$, where $A$ is the uniform closure of $H^\infty|_{\Pi_{c;\delta}^n}$, 
for each $g$ satisfying \eqref{e3.0} there exists a function $h\in (H^\infty)^{-1}$ (depending on $g$ and the data) such that
\begin{equation}\label{e3.4}
\bigl\|(h- g)|_{\Pi_{c;\delta}^n}\bigr\|_{H^\infty(\Pi_{c;\delta}^n)}\le\varepsilon .
\end{equation}
By $\mathcal H_{g,\mathcal D}$ we denote the class of such functions $h$. 

We have to prove that
\begin{equation}\label{e3.5}
C=C(c,M, n,\delta,\varepsilon):=\sup_{F,g}\inf_{h\in \mathcal H_{g,\mathcal D}} \max\bigl\{\|h\|_{H^\infty}, \|h^{-1}\|_{H^\infty}\bigr\}
\end{equation}
is finite. 

To this end, let $\{\Pi_c^n[F_i]\}_{i\in\N}\subset D_{-}$, $F_i=\{f_{ki}\}_{k=1}^n\subset A(\Di)$, and $\{g_i\in (H^\infty(\Pi_c^n[F_i]))^{-1}\}_{i\in\N}$ be sequences of polyhedra and functions satisfying assumptions of the theorem such that
\begin{equation}\label{e3.5a}
C=\lim_{i\rightarrow\infty}\inf_{h\in \mathcal H_{g_i,\mathcal D_i}} \max\bigl\{\|h\|_{H^\infty}, \|h^{-1}\|_{H^\infty}\bigr\};
\end{equation}
here $\mathcal D_i:=\{F_i, c, M, n,\delta,\varepsilon\}$.

\noindent We define functions $f_{k}\in H^\infty(\Di\times\N)$, $1\le k\le n$, by the formulas
\begin{equation}\label{e3.6}
f_k(z,i):=f_{ki}(z),\quad (z,i)\in\Di\times\N.
\end{equation}
Next, we set $F=\{f_k\}_{k=1}^n$ and define polyhedra  $\Pi_{\nu}[F]\subset\Di\times\N$,
\begin{equation}\label{e3.7}
\Pi_{\nu}[F]:=\{x\in\Di\times\N\, :\, \max_{1\le k\le n}|f_k(x)|<\nu \}\, \bigl(=\{(z,i)\in\Di\times\N\, :\, z\in \Pi_{c;\nu}^n[F_i]\}\bigr).
\end{equation}
Finally, we define $g\in H^\infty(\Pi_{1}[F])$,
\begin{equation}\label{e3.8}
g(z,i):=g_i(z),\quad (z,i)\in \Pi_{1}[F].
\end{equation}

Let $\tau:\Z\rightarrow\N$ be a bijection. Consider the biholomorphic map
\begin{equation}\label{e3.9}
b:=({\rm id}_{\Di}\times\tau)\circ s: U\rightarrow\Di\times\N, 
\end{equation}
see \eqref{equ3.1a}, \eqref{equ3.2} in Section~2.2.
We pull back functions $f_k$ and $g$ to $U$ by $b$.
Hence, $b^*f_k:=f_k\circ b\in H^\infty(U)$ and $b^*g\in H^\infty(b^{-1}(\Pi_1[F]))$. By definition, $b^{-1}(\Pi_\nu [F])$ are polyhedra in $U$ determined by functions $b^*f_k$, $1\le k\le n$, cf. \eqref{e3.7}.
Let $\widehat{b^*f_k}\in\mathcal O(\widetilde U)$ be the extension of $b^*f_k$ to $\widetilde U=\hat r^{-1}(\Di)\subset\mathfrak M(H^\infty)$, see \eqref{equ3.4}. We define
open polyhedra $\hat{\Pi}_\nu\subset\widetilde U$ by the formulas
\begin{equation}\label{e3.10}
\hat{\Pi}_{\nu}:=\bigl\{x\in\widetilde U\, :\, \max_{1\le k\le n}|\widehat{b^*f_k}(x)|<\nu \bigr\}.
\end{equation}
Then the open polyhedron $b^{-1}(\Pi_\nu[F])$ is dense in $\hat{\Pi}_{\nu}$. In particular,
by \cite[Thm.\,3.2]{S1}, $b^*g$ is extended to a function $\widehat{b^*g}\in \mathcal O(\hat{\Pi}_1)$. By assumptions of the theorem  $(\widehat{b^*g})^{-1}\in \mathcal O(\hat{\Pi}_1)$ as well. 

Further, since each  $\Pi_c^n[F_i]\subset D_{-}$, the open polyhedron $\hat{\Pi}_{1}$ lies in the $\mathcal O(\mathfrak M(H^\infty))$-convex compact set $N$, see  \eqref{eq2.6}. Moreover, every open polyhedron $\hat{\Pi}_{\nu}$, $0<\nu<1$, is a relatively compact subset of $\hat{\Pi}_{1}$. Hence, due to Lemma \ref{lemma2.2} the $\mathcal O(\widetilde U)$-convex set $K_\delta:=\cap_{\nu>\delta}\hat{\Pi}_{\nu}$ (-- the $\mathcal O(\widetilde U)$-convex hull of $\hat{\Pi}_{\delta}$) is $\mathcal O(\mathfrak M(H^\infty))$-convex.
Since function $\widehat{b^*g}$ is defined on an open neighbourhood of $K_\delta$, \cite[Cor.\,2.6]{S1} implies that $\widehat{b^*g}$ can be uniformly approximated on $K_\delta$ by invertible functions from $\mathcal O(\mathfrak M(H^\infty))$. Thus there is an invertible function $H\in \mathcal O(\mathfrak M(H^\infty))$ such that
\begin{equation}\label{e3.11}
\sup_{x\in \hat{\Pi}_{\delta}}\bigl|H(x)-\widehat{b^*g}(x)\bigr|\le\varepsilon .
\end{equation} 
We set
\begin{equation}\label{e3.12}
h_{i}(z):= (H|_{U}\circ b^{-1})(z,i),\quad (z,i)\in\Di\times\N.
\end{equation}
Then $h_i\in (H^\infty)^{-1}$ and 
\[
\| (h_i-g_i)|_{\Pi_{c;\delta}^n[F_i]}\|_{H^\infty}\le\varepsilon\quad {\rm for\ all}\quad i\in\N.
\]
Moreover,
\[
\sup_{i\in\N}\max\bigl\{\|h_i\|_{H^\infty}, \|h_i^{-1}\|_{H^\infty}\bigr\}\le \max\bigl\{ \|H|_{\Di}\|_{H^\infty},\|H^{-1}|_{\Di}\|_{H^\infty}\bigr\}.
\]
Hence,
\[
C\le \max\bigl\{ \|H|_{\Di}\|_{H^\infty},\|H^{-1}|_{\Di}\|_{H^\infty}\bigr\}
\]
as well, see \eqref{e3.5a}.

This completes the proof of the theorem.
\end{proof}
\begin{proof}[{\bf 3.2.} Proof of Theorem \ref{te1}]
The corona theorem for $H^\infty(\Di\times\N)$ follows from Carleson estimates for solutions of the corona problem for $H^\infty$, see, e.g., \cite[Lm.\,1]{Be}. Hence, $\Di\times\N$ is embedded into $\mathfrak M(H^\infty(\Di\times\N))$ as an open dense subset in the Gelfand topology. 

A compact subset $K\subset \mathfrak M(H^\infty(\Di\times\N))$ is called {\em holomorphically convex} if for each $x\not\in K$ there is $f\in H^\infty(\Di\times\N)$ such that
\begin{equation}\label{e3.14}
\max_{K}|\hat f|< |\hat f(x)|;
\end{equation}
here $\, \hat{\,}:H^\infty(\Di\times\N)\rightarrow C(\mathfrak M(H^\infty(\Di\times\N))$ is the Gelfand transform.

Let $\varphi: H^\infty(\Di\times\N)\rightarrow A$ be a dense morphism of complex unital Banach algebras.  Then the transpose map $\varphi^*$ establishes a homeomorphism from $\mathfrak M(A)$ onto a holomorphically convex subset of $\mathfrak M(H^\infty(\Di\times\N))$. Without loss of generality we identify $\mathfrak M(A)$ with its image under $\varphi^*$.  In order to prove the theorem we have to show that for every $g\in H^\infty(\Di\times\N)$ such that $\hat g$ is zero free on $\mathfrak M(A)$ there is a sequence $\{h_i\}_{i\in\N}\subset (H^\infty(\Di\times\N))^{-1}$  such that $\{\hat h_i|_{\mathfrak M(A)}\}$ converges uniformly to $\hat g|_{\mathfrak M(A)}$, see, e.g., \cite[page\, 262]{S1} and references therein.

To this end, let $O\subset \mathfrak M(H^\infty(\Di\times\N))$ be an open neighbourhood of $\mathfrak M(A)$ such that $\hat{g}$ is yet zero free on its closure $\bar{O}$. Since  $\mathfrak M(A)$ is holomorphically convex, there exist functions $f_{k}\in H^\infty(\Di\times\N)$, $1\le k\le n$, $n\in\N$, and $\delta\in (0,1)$ such that $\mathfrak M(A)\subset
\Pi_{\delta}\subset\Pi_1\subset O$, where
\begin{equation}\label{e3.15}
\Pi_\nu:=\bigl\{x\in \mathfrak M(H^\infty(\Di\times\N))\, :\, \max_{1\le k\le n}|\hat f_k(x)|<\nu\bigr\} ,
\end{equation}
see \cite[Lm.\,5.1]{Br1} for a similar argument.

We set 
\[
g_i(z):=g(i,z),\quad  f_{ki}(z):=f_k(i,z),\quad (i,z)\in\Di\times\N.
\]
If $F_i=\{f_{ki}\}_{1\le k\le n}\subset H^\infty$, then 
\[
\Pi_{\nu}\cap (\Di\times\{i\})=\Pi_{c;\nu}^n[F_i],
\]
where $c:=\max_{1\le k\le n}\|f_k\|_{H^\infty(\Di\times\N)}$, see definition \eqref{e3.1a}.

Moreover,
\[
\max\left\{\|g_i\|_{H^\infty(\Pi_{c;\delta}^n[F_i])},\|g_i^{-1}\|_{H^\infty(\Pi_{c;\delta}^n[F_i])} \right\}\le M:=\max_{\bar{O}}\{|\hat g|,|\hat g)|^{-1}\}.
\]
Thus we can apply Theorem \ref{lem3.1} to functions $g_i$. According to this theorem, there exist
sequences $\{h_{ij}\}_{j\in\N}\in (H^\infty)^{-1}$, $i\in\N$, such that  
$\{h_{ij}|_{\Pi_{c;\delta}^n[F_i]}\}_{j\in\N}$ converges uniformly to $g_i|_{\Pi_{c;\delta}^n[F_i]}$
and for each $j\in\N$
\[
\sup_{i\in\N}\max\bigl\{\|h_{ij}\|_{H^\infty},\|h_{ij}^{-1}\|_{H^\infty}\bigr\}<\infty .
\]
We set
\[
h_j(z,i):=h_{ij}(z),\quad (z,i)\in\Di\times\N,\quad j\in\N.
\]
Then each $h_j\in (H^\infty(\Di\times\N))^{-1}$ and $\{h_j|_{\Pi_\delta}\}_{j\in\N}$ converges uniformly to $g|_{\Pi_{\delta}}$.

This completes the proof of the theorem.
\end{proof}
\sect{Proofs of Proposition \ref{prop1.3} and Theorem \ref{te4.1}}
\begin{proof}[{\bf 4.1.} Proof of Proposition \ref{prop1.3}]
The proof follows the lines of the proof of the Theorem in \cite[Sec.\,2.7]{Ta}, where instead of the Ramspott theorem used in the proof of Proposition 2.6 of that paper one uses a similar result for maps between Stein and Oka manifolds, see \cite[Sec.\,5.4]{F}.
\end{proof}
\noindent {\bf 4.2.} We recall some definitions and results used in the proof of Theorem \ref{te4.1} (for details see, e.g., \cite{Ta}).

For a complex unital Banach algebra $A$ we denote by $\ \hat{\, }:A\hookrightarrow A^{**}$ the natural embedding. Each functional $\hat{a}$, $a\in A$, determines an element of $\mathcal O(\mathfrak M(A))$, see Section~1.2.  We denote by $P(A^*)$ the algebra (under pointwise operations)
generated by the collection of linear functions $\hat{a}$, $a\in A$. 
Then the image of the natural monomorphism $P(A^*)\rightarrow\mathcal O(\mathfrak M(A))$ is a dense subset of $\mathcal O(\mathfrak M(A))$ (see, e.g., \cite[Sect.\,2.3]{Ta}). 

Each $p\in P(A^*)$ can be presented in the form
$p=\sum_{i, \alpha_i}c_{\alpha_i}\hat{a}_i^{\alpha_i}$, where all $c_{\alpha_i}\in\Co$ and the family $\{\hat{a}_i\}\subset A^{**}$ is linearly independent. Elements of this family are referred to as {\em variables} of $p$ so we write $p=p(\hat{a}_1,\hat{a}_2,\dots )$.
If $\{p_j\}\subset P(A^*)$ is a finite family of polynomials, then its set of variables is defined as the maximal linear independent subset of the set of variables of all $p_j$.

Let $K\subset A^*$ be a weak$^*$ compact subset. 
The polynomially convex hull of $K$ is the set
\begin{equation}\label{hull}
\widehat{K}:=\bigl\{y\in A^*\, :\, |p(y)|\le\max_{K}|p|\quad  \forall\, p\in P(A^*)\bigr\}.
\end{equation}
It is a weak$^*$ compact subset of $A^*$ as well.
Such set $K$ is said to be {\em polynomially convex} if $\widehat{K}=K$. For instance, $\mathfrak M(A)\subset A^*$ is a polynomially convex set.

 Further, a {\em polynomial polyhedron} is a set of the form 
 \begin{equation}\label{polyh}
 \Pi =\{y\in A^*\, :\,  \max_{1\le i\le n}|p_i(y)| < 1\}\quad {\rm for}\quad p_1,\dots, p_n\in P(A^*).
 \end{equation}
For each weak$^*$ open neighbourhood $U$ of a polynomially convex set $K$ there is an open polynomial polyhedron $\Pi$ such that $K\subset\Pi\subset U$ (see, e.g., \cite[Sec.\,2.2]{Ta}).
 
Recall that $T_A^n:\mathcal O(\mathfrak M(A),\Co^n)\rightarrow A^n$ denotes the holomorphic functional calculus and $A_\mathcal M:=T_A^n(\mathcal O(\mathfrak M(A),\mathcal M))$ for a complex submanifold $\mathcal M\subset \Co^n$, see \eqref{eq1.4}.
 
 Let $r: A_1\rightarrow A_2$ be a dense image morphism of  commutative complex unital Banach algebras. Then the transpose map
$r^*: A_2^*\rightarrow A_1^*$ is a bounded linear injection.
Without loss of generality we identify $A_2^*$ with its image under $r^*$ so that $r$ is identified with the restriction map $\hat{a}\mapsto\hat{a}|_{A_2^*}$, $a\in A_1$.
The density of the image of $r$ implies that $\mathfrak M(A_2)$ is a polynomially convex subset of $\mathfrak M(A_1)$. 
Also, the restriction map to $A_{2}^*$ maps $\mathcal O(\mathfrak M(A_1),\Co^n)$ in $\mathcal O(\mathfrak M(A_2),\Co^n)$.

\begin{Prop}\label{prop4.0}
We have for all $f\in\mathcal O(\mathfrak M(A_1),\Co^n)$,
\[
(r^n\circ T_{A_1}^n)(f)=T_{A_2}^n(f|_{A_2^*}).
\]
Here $r^n:=(r,\dots ,r):A_1^n\rightarrow A_2^n$.
\end{Prop}
\begin{proof}
It suffices to prove the result for $n=1$ and $f\in P(A_1^*)$ (because $P(A_1^*)$ forms a dense subset of $\mathcal O(\mathfrak M(A_1))$). So assume that $f$ is a polynomial in variables $\hat{v}_{1},\dots, \hat{v}_k$ for some $v_1,\dots, v_k\in A_1$. Then
\[
(r\circ T_{A_1}^1)(f(\hat{v}_{1},\dots, \hat{v}_k))=r(f(v_1,\dots, v_k))=
f(r(v_1),\dots, r(v_k))=T_{A_2}^1(\hat{v}_{1}|_{A_2^*},\dots, \hat{v}_k|_{A_2^*}),
\]
as required.
\end{proof}
{\bf 4.3.} 
Let $\mathcal M\subset\Co^n$ be a complex submanifold and an Oka manifold.  Proposition \ref{prop4.0} implies that
$r^n$ maps $(A_1)_{\mathcal M}$ in $(A_2)_{\mathcal M}$. Let $\overline{r^n((A_1)_{\mathcal M})}$ be the closure of $r^n((A_1)_{\mathcal M})$ in $(A_2)_M$. Theorem \ref{te4.1}\,(1) asserts that
\[
\overline{r^n((A_1)_\mathcal M)}=\bigl\{T_{A_2}^n(f),\, f\in\mathcal O(\mathfrak M(A_2),\mathcal M)\, :\, \exists\, \tilde f\in C(\mathfrak M(A_1),\mathcal M),\  \tilde f|_{\mathfrak M(A_2)}=f|_{\mathfrak M(A_2)}\bigr\}.
\]
The most difficult part of the proof  is to establish that \medskip 

\noindent $(\circ)$ {\em If  $f\in\mathcal O(\mathfrak M(A_2),\mathcal M)$ is such that $f|_{\mathfrak M(A_2)}$ is extended to a map $\tilde f\in C(\mathfrak M(A_1),\mathcal M)$, then} $T_{A_2}^n(f)\in \overline{r^n((A_1)_\mathcal M)}$.\medskip 

In this part, we show that it suffices to prove this statement for some special $f$.\smallskip

Without loss of generality we may assume that $f\in\mathcal O(\Pi,\mathcal M)$ for an open polynomial polyhedron $\Pi\subset\mathfrak M(A_2)$ determined by polynomials from $P(A_2^*)$ in variables $\hat{a}_1,\dots, \hat{a}_m$. Then $\Pi=\pi_m^{-1}(V)$, where 
$\pi_m=(\hat{a}_1,\dots, \hat{a}_m):A_2^*\rightarrow\Co^m$ is surjective and $V$ is an ordinary open polynomial polyhedron in $\Co^m$ containing $\pi_m(\mathfrak M(A_2))$.  Shrinking $\Pi$, if necessary, we may assume that $f$ is a bounded function and there is a bounded function $g\in \mathcal O(V,\mathcal M)$ such that $\pi_m^*(g)=f$, see, e.g., Proposition in \cite[page~159]{Ta}.

Further, by the density of the image of $r$, there is a sequence $\{(b_{1j},\dots,b_{mj})\}_{j\in\N}\subset A_1^m$ such that the sequence $\{(r(b_{1j}),\dots,r(b_{mj}))\}_{j\in\N}\subset A_2^m$ converges to $(a_1,\dots,a_m)$. Then the sequence of maps $\{\pi_{mj}|_{A_2^*}\}_{j\in\N}$, where
$\pi_{mj}:=(\hat{b}_{1j},\dots,\hat{b}_{mj}):A_1^*\rightarrow\Co^m$,  converges to $\pi_m$ in the norm topology of the space of bounded linear maps $A_2^*\rightarrow\Co^n$.
 Hence, there is $j_0\in\N$ such that $\pi_{mj}(\mathfrak M(A_2))\subset V$  and maps $\pi_{mj}|_{A_2^*}$ are surjective for all $j\ge j_0$. In particular, $g(\pi_{mj}(\mathfrak M(A_2)))\subset \mathcal M$ for such $j$. 
Let us consider functions $\pi_{mj}^*(g)\in\mathcal O(\Pi_j,\mathcal M)$, $j\ge j_0$; here $\Pi_j:=\pi_{mj}^{-1}(V)\subset A_1^*$ are open polynomial polyhedra containing $\mathfrak M(A_2)$. Since the map $r:A_2^*\hookrightarrow A_1^*$ is weak$^*$ continuous, $\Pi_j\cap A_2^*$ are weak$^*$ open neighbourhoods of $\mathfrak M(A_2)$ in $A_2^*$. Hence, $\pi_{mj}^*(g)|_{A_2^*}\in\mathcal O(\mathfrak M(A_2),\mathcal M)$ and $T_{A_2}^n(\pi_{mj}^*(g)|_{A_2^*})\in (A_2)_\mathcal M$.

\begin{Lm}\label{lemma4.2}
The sequence $\{T_{A_2}^n(\pi_{mj}^*(g)|_{A_2^*})\}_{j\ge j_0}$ converges to $T_{A_2}^n(f)$.
\end{Lm}
\begin{proof}
For the basic results of the theory of Stein manifolds see, e.g., \cite{GR}.

Since $V\subset\Co^m$ is a polynomial polyhedron, there is a sequence of holomorphic polynomials $\{g_i\}_{i\in\N}$ on $\Co^m$ which converges to $g$ uniformly on compact subsets of $V$. By the definition $\{\pi_{mj}^*(g_i)|_{A_2^*}\}_{i\in\N}\subset P(A_2^*)$ is a sequence of polynomials in variables $\hat{b}_{1j}|_{A_2^*},\dots, \hat{b}_{mj}|_{A_2^*}$, where $\hat{b}_{kj}|_{A_2^*}=\widehat{r(b_{kj})}$, $1\le k\le m$, converging to $\pi_{mj}^*(g_i)|_{A_2^*}$ uniformly on preimages under $\pi_{mj}|_{A_2^*}$ of compact subsets of $V$ with interiors containing $\pi_{mj}(\mathfrak M(A_2))$. This implies that
\[
T_{A_2}^n(\pi_{mj}^*(g)|_{A_2^*})=\lim_{i\rightarrow\infty}T_{A_2}^n(\pi_{mj}^*(g_i)|_{A_2^*})=
\lim_{i\rightarrow\infty} g_i(r(b_{1j}),\dots, r(b_{mj})).
\]
Thus given $\varepsilon>0$ there is $i_0\in\N$ such that for each $i\ge i_0$,
\begin{equation}\label{eq4.1}
\| T_{A_2}^n(\pi_{mj}^*(g)|_{A_2^*})- g_i(r(b_{1j}),\dots, r(b_{mj}))\|_{A_2^n}\le\frac{\varepsilon}{3}.
\end{equation}

Next, since $\lim_{j\rightarrow\infty}r(b_{kj})=a_k$, $1\le k\le m$,  and $g_i$ is a polynomial,  $\newline \lim_{j\rightarrow\infty}g_i(r(b_{1j}),\dots, r(b_{mj}))=g_i(a_1,\dots, a_m)$. Hence, there is $j_1\in\N$, $j_1\ge j_0$, such that for each $j\ge j_1$,
\begin{equation}\label{eq4.2}
\|g_i(r(b_{1j}),\dots, r(b_{mj}))-g_i(a_1,\dots, a_m)\|_{A_2^n}\le\frac\varepsilon 3 .
\end{equation}

Finally,  
\[
\lim_{i\rightarrow\infty}g_i(a_1,\dots, a_m)=\lim_{i\rightarrow\infty}T_{A_2}^n ( \pi_m^*(g_i))=
T_{A_2}^n\left(\pi_m^*\left(\lim_{i\rightarrow\infty}g_i\right)\right)=T_{A_2}^n(f).
\]
Hence, there is $i_1\ge i_0$ such that for each $i\ge i_1$
\begin{equation}\label{eq4.3}
\| g_i(a_1,\dots, a_m)-T_{A_2}^n(f) \|_{A_2^n}\le\frac\varepsilon 3 .
\end{equation}
Adding inequalities \eqref{eq4.1}--\eqref{eq4.3} we obtain for each $j\ge j_1$,
\[
\|T_{A_2}^n(\pi_{mj}^*(g)|_{A_2^*})-T_{A_2}^n(f)\|_{A_2^n}\le\varepsilon .
\]
This proves the required statement.
\end{proof}
Next, we prove the following result.
\begin{Lm}\label{lemma4.3}
If $f|_{\mathfrak M(A_2)}$ is  extended to a map $\tilde f\in C(\mathfrak M(A_1),\mathcal M)$, then for all sufficiently large $j$ each map $\pi_{mj}^*(g)|_{\mathfrak M(A_2)}$ is extended to a map from $C(\mathfrak M(A_1),\mathcal M)$.
\end{Lm}
\begin{proof}
By the definition of a complex submanifold of $\Co^n$, for each $z\in\mathcal M$ there is an open complex Euclidean ball $B_z\subset\Co^n$ centered at $z$ such that $B_z\cap\mathcal M$ is a closed complex submanifold of $B_z$. This implies that $\mathcal M$ is a {\em closed} complex submanifold of the open set $\underset{z\in\mathcal M}\cup B_z\subset\Co^n$. Indeed, we have 
\[
\left(\underset{z\in\mathcal M}\cup B_z\right)\setminus\mathcal M=\left(\underset{z\in\mathcal M}\cup B_z\setminus (B_z\cap \mathcal M)\right).
\]
Then since each $B_z\setminus (B_z\cap \mathcal M)$ is an open subset of $B_z$, the complement of $\mathcal M$ in $\underset{z\in\mathcal M}\cup B_z$ is open as claimed.

We obtain from here that since $\mathcal M$ is an absolute neighbourhood retract,  
there is an open neighbourhood $U\subset \underset{z\in\mathcal M}\cup B_z$ of $\mathcal M$ and a continuous retraction $R: U\rightarrow \mathcal M$. 

Further, by our definition each $\pi_{mj}^*(g)$, $j\ge j_0$, maps the weak$^*$ open neighbourhood $\Pi_j\subset A_1^*$ of $\mathfrak M(A_2)$ into $\mathcal M$. Also, Lemma \ref{lemma4.2} implies that the sequence $\{\pi_{mj}^*(g)|_{A_2^*}\}_{j\ge j_0}$ converges to $f$ uniformly on $\mathfrak M(A_2)$. Hence, given $\varepsilon>0$ there is $j_1\in\N$, $j_1\ge j_0$, and for each $j\ge j_1$ there is an open in the Gelfand topology of $\mathfrak M(A_1)$ neighbourhood $O_{j,\varepsilon}\subset \mathfrak M(A_1)$ of $\mathfrak M(A_2)$ contained in $\Pi_j$ such that 
\begin{equation}\label{eq4.4}
\|\pi_{mj}^*(g)(y)-f(y)\|_{\Co^n}\le\varepsilon\quad {\rm for\ all}\quad y\in O_{j,\varepsilon}.
\end{equation}
Let $\rho_1\in C(O_{j,\varepsilon})$, $\rho_2\in C(\mathfrak M(A_1)\setminus\mathfrak M(A_2))$ be a continuous partition of unity subordinate to the open cover $\{O_{j,\varepsilon},\mathfrak M(A_1)\setminus\mathfrak M(A_2)\}$ of the compact space $\mathfrak M(A_1)$. Then \eqref{eq4.4} implies
\begin{equation}\label{eq4.5}
\sup_{y\in\mathfrak M(A_1)}\| \rho_1(y) \pi_{mj}^*(g)(y)+\rho_2(y) f(y)-f(y)   \|_{\Co^n}\le\varepsilon .
\end{equation}
Since $f(\mathfrak M(A_1))\subset \mathcal M$ is a compact subset, \eqref{eq4.5} implies that for a sufficiently small  $\varepsilon$,  $(\rho_1\pi_{mj}^*(g)+\rho_2 f)(y)\in U$ for all $y\in\mathfrak M(A_1)$. We define
\[
h_j:=R\circ (\rho_1\pi_{mj}^*(g)+\rho_2 f).
\]
Then $h_j\in C(\mathfrak M(A_1),\mathcal M)$ and since $\rho_2=0$ on $\mathfrak M(A_2)$, 
\[
h_j|_{\mathfrak M(A_2)}=R\circ (\pi_{mj}^*(g)|_{\mathfrak M(A_2)})=\pi_{mj}^*(g)|_{\mathfrak M(A_2)}.
\]
This completes the proof of the lemma.
\end{proof}

The above arguments  show that it suffices to prove statement $(\circ)$ for elements $\pi_{mj}^*(g)|_{A_2^*}$ for all sufficiently large $j$ in place of $f$.

\begin{proof}[{\bf 4.4.} Proof of Theorem \ref{te4.1}(1)]
Let  $f\in\mathcal O(\mathfrak M(A_2),\mathcal M)$ be such that $\tilde f|_{\mathfrak M(A_2)}=f|_{\mathfrak M(A_2)}$ for some $\tilde f\in C(\mathfrak M(A_1),\mathcal M)$. We have to prove that $T_{A_2}^n(f)\in \overline{r^n((A_1)_\mathcal M)}$. \smallskip

Due to the reduction presented in the previous section, it suffices to prove this result for a bounded function $f\in\mathcal O(\Pi, \mathcal M)$, where $\Pi$ is an {\em open polynomial polyhedron in} $A_1^*$ containing $\mathfrak M(A_2)$. Shrinking $\Pi$, if necessary, we may assume that $f(\Pi)$ is a relatively compact subset of $\mathcal M$. By the Stone-Weierstrass theorem, the uniform algebra generated by restrictions of $P(A_1^*)$ and its complex conjugate $\bar{P}(A_1^*)$ to $\mathfrak M(A_1)$ coincides with $C(\mathfrak M(A_1))$. Hence, given $\varepsilon>0$ there is a (nonholomorphic) polynomial
map $p_\varepsilon: A_1^*\rightarrow\Co^n$ such that
\begin{equation}\label{equ4.6}
\sup_{x\in\mathfrak M(A_1)}\|p_\varepsilon(x)-\tilde f(x)\|_{\Co^n}<\varepsilon .
\end{equation}
Since by our hypothesis the closure $\overline{f(\Pi)}$ of $f(\Pi)$ is a compact subset of $\mathcal M$, we can choose $\varepsilon>0$ so small that the open $\varepsilon$-neighbourhood of $\overline{f(\Pi)}$,
\[
\overline{f(\Pi)}_\varepsilon:=\{z\in\Co^n\, :\, \exists\, x\in\Pi,\ \|f(x)-z\|_{\Co^n}<\varepsilon\},
\]
is contained in the open neighbourhood $U$ of Lemma \ref{lemma4.3}.  
Since the map $p_\varepsilon|_{\Pi}-f\in C(\Pi,\Co^n)$ is weak$^*$ continuous and $\Pi$ is a weak$^*$ open subset of $A_1^*$, the set $(p_\varepsilon|_{\Pi}-f)^{-1}(B_\varepsilon)$ is a weak$^*$ open subset of $A_1^*$ containing $\mathfrak M(A_2)$; here $B_\varepsilon:=\{z\in\Co^n\, :\, \|z\|_{\Co^n}<\varepsilon\}$. In particular, $\Pi\cap \bigl((p_\varepsilon|_{\Pi}-f)^{-1}(B_\varepsilon)\bigr)$ is a weak$^*$ open neighbourhood of $\mathfrak M(A_2)$ and so it contains an open polynomial polyhedron $\widetilde \Pi_\epsilon$ with the same property. Then we have 
\begin{equation}\label{equ4.6a}
\lambda p_\varepsilon(x)+(1-\lambda)f(x)\in U\quad {\rm for\ all}\quad x\in \widetilde \Pi_\epsilon,\ \lambda\in [0,1].
\end{equation}

Further, due to \eqref{equ4.6}, $p_\varepsilon^{-1}(U)$ is a weak$^*$ open neighbourhood of $\mathfrak M(A_1)$. Then there is an open polynomial polyhedron $\Pi_\varepsilon\subset p_\varepsilon^{-1}(U)$ containing $\mathfrak M(A_1)$ such that $p_\varepsilon\in C(\Pi_\varepsilon,\Co^n)$ is a bounded map. In turn, $\Pi_\varepsilon':=\widetilde\Pi_\varepsilon\cap \Pi_\varepsilon$ is an open polynomial polyhedron containing $\mathfrak M(A_2)$.   Suppose that $\Pi_\varepsilon'$ is determined by polynomials in variables $\hat{v}_{1},\dots, \hat{v}_l$ for some $v_1,\dots, v_l\in A_1$ (so that the variables of polynomials defining $\Pi_\varepsilon$ are linear combinations of  these  variables).  Consider the bounded linear surjective map $\pi_l=(\hat v_1,\dots, \hat v_l): A_1^*\rightarrow\Co^{l}$. Then there are ordinary open polynomial polyhedra $Q_\varepsilon'\subset Q_\varepsilon\subset\Co^l$ such that $\pi_l^{-1}(Q_\varepsilon')=\Pi_\varepsilon'$ and
$\pi_l^{-1}(Q_\varepsilon)=\Pi_\varepsilon$. Also, there are bounded maps $q_\varepsilon\in C(Q_\varepsilon' ,\Co^n)$  and  $g_\varepsilon\in\mathcal O(Q_\varepsilon', \Co^n)$ such that 
\[
\pi_l^*(q_\varepsilon)=p_\varepsilon\quad {\rm and}\quad \pi_l^*(g_\varepsilon)=f|_{\Pi_\varepsilon'}.
\]
Note that $K:=\pi_l(\mathfrak M(A_2))$ is a compact subset of $Q_\varepsilon'$. Hence the polynomially  convex hull $\widehat{K}$ of $K$ is a compact subset of $Q_\varepsilon'$ as well.  We choose a compact polynomially convex subset $S\Subset Q_{\varepsilon}'$ whose interior contains $\widehat{K}$.
Let $\rho_1\in C(Q_\varepsilon')$ and $\rho_2\in C(Q_\varepsilon\setminus S )$ be a continuous partition of unity subordinate to the open cover $\{Q_\varepsilon',Q_\varepsilon\setminus S\}$ of $Q_\varepsilon$ such that $\rho_1=1$ on an open neighbourhood of $S$. Then by our choice of $\varepsilon$, see \eqref{equ4.6a},
 the map
\[
h_\varepsilon:=R\circ (\rho_1 g_\varepsilon+\rho_2q_\varepsilon)\in C(Q_\varepsilon,\mathcal M)
\]
and coincides with $g_\varepsilon$ on an open neighbourhood of $S$; in particular, $h_\varepsilon$ is holomorphic there.

Recall that each Oka manifold $Y$ satisfies the
 {\em basic Oka property with approximation}. This means that every continuous map $f_0:X\rightarrow Y$ from a Stein space $X$ that is holomorphic on (a neighbourhood of) a compact $\mathcal O(X)$-convex subset $C \subset X$ can be deformed to a holomorphic map $f_1:X\rightarrow Y$ by a homotopy of maps that are holomorphic near $C$ and arbitrary close to $f_0$ on $C$, \cite[Sect.\,5.15]{F}.

We apply this property to our case with $X:=Q_\varepsilon$, $Y:=\mathcal M$, $C:=S$ and $f_0:=h_\varepsilon$. Then we obtain that there is a sequence of holomorphic maps $\{F_i\}_{i\in\N}\subset\mathcal O(Q_\varepsilon,\mathcal M)$ such that
\[
\lim_{i\rightarrow\infty}\max_{z\in S}\|F_i(z)-g_\varepsilon(z)\|_{\Co^n}=0.
\]
Hence, the sequence $\{\pi_l^*(F_i)\}_{i\in\N}\in \mathcal O(\Pi_\varepsilon,\mathcal M)$ converges uniformly to $f$ on the open neighbourhood $\pi_l^{-1}(\mathring{S})$ of $\mathfrak M(A_2)$. (Here $\mathring{S}$ is the interior of $S$.) We set
\[
c_i:=T_{A_1}^n(F_i),\quad i\in\N.
\]
Then each $c_i\in (A_1)_\mathcal M$ and 
\[
\lim_{i\rightarrow\infty} r^n(c_i)=\lim_{i\rightarrow\infty}T_{A_2}^n(F_i|_{A_2^*})=T_{A_2}^n\left(\lim_{i\rightarrow\infty}F_i|_{A_2^*}\right)=T_{A_2}^n(f).
\]
This shows that $T_{A_2}^n(f)\in \overline{r^n((A_1)_\mathcal M)}$ which completes the proof of the first part of Theorem \ref{te4.1}\,(1).

Conversely, suppose that $T_{A_2}^n(f)$, $f\in\mathcal O(\mathfrak M(A_2),\Co^n)$, belongs to $\overline{r^n((A_1)_\mathcal M)}$. We have to show that there is $\tilde f\in C(\mathfrak M(A_1),\Co^n)$ such that
$\tilde f|_{\mathfrak M(A_2)}=f$.

Indeed, since $T_{A_2}^n(f)\in \overline{r^n((A_1)_\mathcal M)}$, there is a sequence $\{f_i\}_{i\in\N}\subset \mathcal O(\mathfrak M(A_1),\Co^n)$ such that the sequence $\{r^n(T_{A_1}^n(f_i))\}_{i\in\N}$ converges to $T_{A_2}^n(f)$. This implies that $\{f_i\}_{i\in\N}$ converges to $f$ uniformly on $\mathfrak M(A_2)$. Since $\mathcal M$ is an absolute neighbourhood retract, $f$ can be extended to a map $f'\in C(O,\mathcal M)$, where $O\subset\mathfrak M(A_1)$ is an open (in the Gelfand topology of $\mathfrak M(A_1)$) neighbourhood of $\mathfrak M(A_2)$. 
Hence, given $\varepsilon>0$ there is $i_0\in\N$ and  an open neighbourhood $O_{i_0,\varepsilon}\subset \mathfrak M(A_1)$ of $\mathfrak M(A_2)$ contained in $O$
such that 
\begin{equation}\label{eq4.8}
\sup_{y\in O_{i_0,\varepsilon}}\|f_{i_0}(y)-f'(y)\|_{\Co^n}\le\varepsilon.
\end{equation}
Let $\rho_1\in C(O_{i_0,\varepsilon})$, $\rho_2\in C(\mathfrak M(A_1)\setminus\mathfrak M(A_2))$ be a continuous partition of unity subordinate to the open cover $\{O_{i_0,\varepsilon},\mathfrak M(A_1)\setminus\mathfrak M(A_2)\}$ of the compact space $\mathfrak M(A_1)$. Then due to \eqref{eq4.8},
\begin{equation}\label{eq4.9}
\sup_{y\in\mathfrak M(A_1)}\|\rho_1(y) f'(y)+ \rho_2(y) f_{i_0}(y)-f_{i_0}(y)   \|_{\Co^n}\le\varepsilon .
\end{equation}
Since $f_{i_0}(\mathfrak M(A_1))\subset \mathcal M$ is a compact set, \eqref{eq4.9} implies that for a sufficiently small  $\varepsilon$,  $\rho_1(y) f'(y)+ \rho_2(y) f_{i_0}(y)\in U$ (see Lemma \ref{lemma4.3}) for all $y\in\mathfrak M(A_1)$.
 We define
\[
\tilde f:=R\circ (\rho_1 f'+\rho_2f_{i_0})\in C(\mathfrak M(A_1),\mathcal M).
\]
Since $\rho_2=0$ on $\mathfrak M(A_2)$, 
\[
\tilde f|_{\mathfrak M(A_2)}=R\circ (f'|_{\mathfrak M(A_2)})=f.
\]
This completes the proof of part (1) of Theorem \ref{te4.1}.
\end{proof}
\begin{proof}[{\bf 4.5.} Proof of Theorem \ref{te4.1}(2)]
We retain notation of the previous sections. Our goal  is to prove that under hypotheses of the theorem
\[
\overline{r^n((A_1)_\mathcal M)}=r^n((A_1)_\mathcal M).
\]

To this end, let $T_{A_2}^n(f)$, $f\in\mathcal O(\mathfrak M(A_2),\mathcal M)$, belong to $\overline{r^n((A_1)_\mathcal M)}$. We have to prove that there exists $\tilde f\in\mathcal O(\mathfrak M(A_1),\mathcal M)$ such that
\begin{equation}\label{equ5.1}
r^n(T_{A_1}^n(\tilde f))= T_{A_2}^n(f).
\end{equation}
As before, we assume that $A_2^*$ is a vector subspace of $A_1^*$.
Since $r$ is a surjective map, the Banach open mapping theorem implies that $A_2^*$ is a {\em weak$^*$ closed} subspace of $A_1^*$. 

Without loss of generality we may assume that $f\in\mathcal O(\Pi,\mathcal M)$ for an open polynomial polyhedron $\Pi\subset\mathfrak M(A_2)$ determined by polynomials from $P(A_2^*)$ in variables $\hat{a}_1,\dots, \hat{a}_m$. Then $\Pi=\pi_m^{-1}(V)$, where 
$\pi_m=(\hat{a}_1,\dots, \hat{a}_m):A_2^*\rightarrow\Co^m$ is surjective and $V$ is an ordinary open polynomial polyhedron in $\Co^m$ containing $\pi_m(\mathfrak M(A_2))$.  Shrinking $\Pi$, if necessary, we may assume that $f$ is a bounded function and there is a bounded function $g\in \mathcal O(V,\mathcal M)$ such that $\pi_m^*(g)=f$.

Further, by surjectivity of $r$, there is an $m$-tuple $(b_{1},\dots,b_{m})\subset A_1^m$ such that  
\[
(r(b_{1}),\dots,r(b_{m}))=(a_1,\dots,a_m).
\]
We set $\tilde\pi_m:=(\hat b_1,\dots,\hat b_m):A_1^*\rightarrow\Co^m$ and consider the function $\tilde\pi_m^*(g)\in \mathcal O(\widetilde\Pi,\mathcal M)$; here $\widetilde\Pi:=\tilde\pi_{m}^{-1}(V)\subset A_1^*$ is an open polynomial polyhedron containing $\mathfrak M(A_2)$. Since, 
$\hat b_i|_{A_2^*}=\hat a_i$, $1\le i\le m$,
\begin{equation}\label{equ5.2}
\tilde\pi_m^*(g)|_{A_2^*}=f.
\end{equation}

Next, since $T_{A_2}^n(f)\in \overline{r^n((A_1)_\mathcal M)}$,
according to Theorem \ref{te4.1}\,(1) there exists a map $f'\in C(\mathfrak M(A_1),\mathcal M)$ such that 
\[
f'|_{\mathfrak M(A_2)}=f.
\]
As in the proof of Theorem \ref{te4.1}\,(1), see Section~4.4, for a sufficiently small $\varepsilon>0$ we can find
 a (nonholomorphic) polynomial
map $p_\varepsilon: A_1^*\rightarrow\Co^n$ and 
open polynomial polyhedra $\widetilde \Pi_\epsilon\subset\widetilde\Pi$ containing $\mathfrak M(A_2)$ and $\hat\Pi_\varepsilon$ containing $\mathfrak M(A_1)$ such that $p_\varepsilon\in C(\hat\Pi_\varepsilon,\Co^n)$ is a bounded map and
\begin{equation}\label{equ5.4}
\begin{array}{c}
\displaystyle \sup_{x\in\mathfrak M(A_1)}\|p_\varepsilon(x)- f'(x)\|_{\Co^n}<\varepsilon,\quad {\rm and}\medskip\\
\displaystyle
\lambda p_\varepsilon(x)+(1-\lambda)\tilde\pi_m^*(g)(x)\in U\quad \forall\, x\in \widetilde \Pi_\epsilon,\, \lambda\in [0,1].
\end{array}
\end{equation}
Without loss of generality we assume that
\begin{equation}\label{equ5.3}
K:=\mathfrak M(A_1)\setminus\widetilde\Pi_\varepsilon\ne\emptyset.
\end{equation}
(For otherwise, $\tilde f:=\tilde\pi_m^*(g)\in\mathcal O(\mathfrak M(A_1),\mathcal M)$ is the required function satisfying \eqref{equ5.1}, see Proposition \ref{prop4.0}.)

By definition, $K$ is a compact subset of $\mathfrak M(A_1)\setminus\mathfrak M(A_2)$; hence,  $K\cap A_2^*=\emptyset$. Then since $A_2^*$ is a weak$^*$ closed subset of $A_1^*$, there is a weak$^*$ open neighbourhood $O\subset A_1^*$ of $K$ such that its weak$^*$ closure $\bar{O}$ does not interest $A_1^*$ as well. In particular,
$\hat\Pi_\varepsilon\cap(O\cup\widetilde\Pi_\varepsilon)\subset A_1^*$ is a weak$^*$ open neighbourhood of $\mathfrak M(A_1)$. Hence, there is an open polynomial polyhedron $\Pi_\varepsilon$ which is contained in this neighbourhood and contains $\mathfrak M(A_1)$.

Continuing as in the proof of Theorem \ref{te4.1}\,(1), we define
an open polynomial polyhedron $\Pi_\varepsilon':=\widetilde\Pi_\varepsilon\cap \Pi_\varepsilon$ in $A_1^*$ containing $\mathfrak M(A_2)$. Then there are a bounded weak$^*$ continuous linear surjective map $\pi_l: A_1^*\rightarrow\Co^l$ and 
ordinary open polynomial polyhedra $Q_\varepsilon'\subset Q_\varepsilon\subset\Co^l$ such that $\pi_l^{-1}(Q_\varepsilon')=\Pi_\varepsilon'$ and $\pi_l^{-1}(Q_\varepsilon)=\Pi_\varepsilon$. Also, there are bounded maps $q_\varepsilon\in C(Q_\varepsilon' ,\Co^n)$  and  $g_\varepsilon\in\mathcal O(Q_\varepsilon', \Co^n)$ such that 
\[
\pi_l^*(q_\varepsilon)=p_\varepsilon\quad {\rm and}\quad \pi_l^*(g_\varepsilon)=\tilde\pi_m^*(g)|_{\Pi_\varepsilon'}.
\]
Let
\begin{equation}\label{equ5.5}
L:=\pi_l(A_2^*)\subset\Co^l.
\end{equation} 
\begin{Lm}
$L$ is a proper linear subspace of $\Co^l$. Moreover,
\begin{equation}\label{equ5.6}
L\cap Q_\varepsilon'=L\cap Q_\varepsilon .
\end{equation}
\end{Lm}
\begin{proof}
By our assumption,
$\Pi_\varepsilon\setminus\Pi_\varepsilon'\ne\emptyset$, see \eqref{equ5.3}, and 
\[
\pi_l^{-1}(Q_\varepsilon\setminus Q_\varepsilon')=\Pi_\varepsilon\setminus\Pi_\varepsilon'.
\] 
Hence, $Q_\varepsilon\setminus Q_\varepsilon'\ne\emptyset$  as well and so it suffices to prove only \eqref{equ5.6}.

Suppose, on the contrary, that $
(L\cap Q_\varepsilon)\setminus(L\cap Q_\varepsilon')\ne\emptyset$. Since
\[
(L\cap Q_\varepsilon)\setminus(L\cap Q_\varepsilon')=L\cap(Q_\varepsilon\setminus Q_\varepsilon'),
\]
there is a point $v\in A_2^*$ such that $\pi_l(v)\in Q_\varepsilon\setminus Q_\varepsilon'$. Then 
\[
v\in \pi_l^{-1}(Q_\varepsilon\setminus Q_\varepsilon')=\Pi_\varepsilon\setminus\Pi_\varepsilon'=\Pi_\varepsilon\setminus\widetilde\Pi_\varepsilon \subset (O\cup\widetilde \Pi_\varepsilon)\setminus \widetilde\Pi_\varepsilon= O\setminus\widetilde\Pi_\varepsilon\subset O.
\]
However, $O\cap A_2^*=\emptyset$ by the choice of $O$. This contradiction proves the lemma.
\end{proof}
The lemma shows that $Z:=L\cap Q_\varepsilon'$ is a closed complex submanifold of $Q_\varepsilon$.

Further,  let $\rho_1\in C(Q_\varepsilon')$ and $\rho_2\in C(Q_\varepsilon\setminus Z )$ be a continuous partition of unity subordinate to the open cover $\{Q_\varepsilon',Q_\varepsilon\setminus Z\}$ of $Q_\varepsilon$ (in particular, $\rho_1=1$ on an open neighbourhood of $Z$). Then by our choice of $\varepsilon$, see \eqref{equ5.4},
 the map
\[
h_\varepsilon:=R\circ (\rho_1 g_\varepsilon+\rho_2q_\varepsilon)\in C(Q_\varepsilon,\mathcal M)
\]
and coincides with $g_\varepsilon$ on an open neighbourhood of $Z$; in particular, $h_\varepsilon$ is holomorphic there.

Recall that each Oka manifold $Y$ satisfies the
 {\em basic Oka property with interpolation}. That is, every continuous map $f_0:X\rightarrow Y$ from a Stein space $X$  holomorphic on an open neighbourhood of a closed complex subvariety $X' \subset X$ can be deformed to a holomorphic map $f_1:X\rightarrow Y$ by a homotopy of maps that is fixed on $X'$, \cite[Sect.\,5.15]{F}.
 
We apply this property to our case with $X:=Q_\varepsilon$, $X'=Z$, $Y:=\mathcal M$ and $f_0:=h_\varepsilon$. Then we obtain a
map $F\in\mathcal O(Q_\varepsilon,\mathcal M)$ such that
\[
F|_Z=g_\varepsilon|_Z.
\]
Hence $\pi_l^*(F)\in \mathcal O(\Pi_\varepsilon,\mathcal M)$ satisfies
\[
\pi_l^*(F)|_{A_2^*}=\pi_l^*(g_\varepsilon)|_{A_2^*}=\tilde\pi_m^*(g)|_{\Pi_\varepsilon'\cap A_2^*}=f|_{\Pi_\varepsilon'\cap A_2^*}.
\]
Since $\Pi_\varepsilon'\cap A_2^*\subset A_2^*$ is an open polyhedron containing $\mathfrak M(A_2)$, the latter and Proposition \ref{prop4.0} imply 
\[
(r^n\circ T_{A_1}^n)(\pi_l^*(F))=T_{A_2}^n(\pi_l^*(F)|_{A_2^*})=T_{A_2}^n(f).
\]
This gives \eqref{equ5.1} with $\tilde f:=\pi_l^*(F)$ and completes the proof of Theorem \ref{te4.1}\,(2).
\end{proof}

\sect{Proof of Theorem \ref{teo1.3}}
\noindent {\bf 5.1.}  First, we prove the following result. 

Let $K\subset \mathfrak M(H^\infty(\Di\times\N))$ be a holomorphically convex set, see \eqref{e3.14}, and let $f\in C(K, \mathcal M)$ be a continuous map into a manifold $\mathcal M$. Let $r:\mathcal M'\rightarrow \mathcal M$ be a finite unbranched covering of $\mathcal M$.
\begin{Prop}\label{prop5.1}
There exists a continuous map $f'\in C(K, \mathcal M')$ such that $f=r\circ f'$.
\end{Prop}
\begin{proof}
Since $\mathcal M$ is an absolute neighbourhood retract, the map $f:K\rightarrow \mathcal M$ admits a continuous extension to a map of an open neighbourhood $O\subset \mathfrak M(H^\infty(\Di\times\N))$  of $K$ into $\mathcal M$. We retain symbol $f$ for the extended map.

Next, we regard $r:\mathcal M'\rightarrow \mathcal M$ as a locally trivial fibre bundle on $\mathcal M$  with a finite fibre $F$. We denote by $\rho : S_d\rightarrow {\rm Aut}(F)$  the natural action of the permutation group $S_d$, $d:={\rm card}\, F$, on fibres of $r$. Let $(U_i)_{i\in I}$ be a cover of $\mathcal M$ by simply connected open sets. Then each $r^{-1}(U_i)$ is homeomorphic to $U_i\times F$. Let
$s_i: U_i\rightarrow \mathcal M'$ be a continuous map (section) such that $r\circ s_i={\rm id}_{U_i}$, $i\in I$. Then there exist  locally constant continuous maps  $g_{ij}: U_i\cap U_j\rightarrow S_d$, $i,j\in I$, such that 
\begin{equation}\label{e4.0}
s_i(x)=\rho(g_{ij}(x))(s_j(x)) \quad {\rm for\ all}\quad  x\in U_i\cap U_j\, (\neq\emptyset).
\end{equation}
We have by the definition,
\[
g_{ij}\cdot g_{jk}\cdot g_{ki}=1 \quad {\rm if}\quad U_i\cap U_j\cap U_k\ne\emptyset\quad {\rm and}\quad g_{ii}=1,\quad g_{ij}=g_{ij}^{-1};
\]
here $1$ is the identity of $S_d$.

Thus $\{g_{ij}\}_{i,j\in I}$ is a locally constant cocycle on the cover $(U_i)_{i\in I}$ with values in $S_d$. In particular, it determines a principal bundle $E\rightarrow \mathcal M$ on $\mathcal M$ with fibre $S_d$ (see, e.g., \cite{Hus} for basic definitions and results of fibre bundles theory). The pullback of $E$ by $f$ is a principal bundle $ f^*E\rightarrow O$ on $O$ with fibre $S_d$ defined on the open cover $(f^{-1}(U_i))_{i\in I}$ by the locally constant cocycle
$\{f^*g_{ij}\}_{i,j\in I}$.  

Let $\Pi\subset O$ be an open polyhedron containing $K$ (existing because $K$ is holomorphically convex), see \eqref{e3.15}. We prove the following result.
\begin{Lm}\label{lem4.2}
The bundle $f^*E$ admits a continuous section over $\Pi$.
\end{Lm}
\begin{proof}
Suppose that 
\[
\Pi=\left\{x\in \mathfrak M(H^\infty(\Di\times\N))\, :\, \max_{1\le k\le n}|\hat f_k(x)|<1\right\} 
\]
for some $f_k\in H^\infty(\Di\times\N)$, $1\le k\le n$. Then 
$\Pi\cap (\Di\times\{i\})$ is an open polyhedron determined by functions$f_{ki}:=f(\cdot,i)\in H^\infty$, $1\le k\le n$, $i\in\N$. By the maximum modulus principle for the subharmonic function $\underset{1\le k\le n}\max |f_{ki}|$, each connected component of $\Pi\cap (\Di\times\{i\})$ is an open simply connected set (hence, it is contractible). In particular, the  bundle $f^*E$ is trivial over each such component and so it is trivial over $\Pi\cap (\Di\times\N)$ which is the disjoint union of these components. Therefore there is a continuous section $g: \Pi\cap (\Di\times\N)$ of $f^*E$. In local coordinates of  $f^*E$ on $(f^{-1}(U_i))_{i\in I}$ this section is given by continuous maps $g_i: f^{-1}(U_{i})\cap (\Di\times\N)\to S_d$, $i\in I$, such that
\begin{equation}\label{e4.1}
g_j(x)=g_i(x)\cdot (f^*g_{ij})(x)\ {\rm for\ all}\, x\in (f^{-1}(U_i)\cap (\Di\times\N))\cap (f^{-1}(U_j)\cap(\Di\times\N)),\  i,j\in I.
\end{equation}

Let $\tau : S_d\subset\N$ be an embedding. Then each $\tau\circ g_i: f^{-1}(U_{i})\cap (\Di\times\N)\to\N$ is a locally constant continuous function with range in $\tau(S_d)$ and, in particular, it is a bounded holomorphic function. Since
each $f^{-1}(U_i)$ is an open subset of  $\mathfrak M(H^\infty(\Di\times\N))$ and $f^{-1}(U_{i})\cap (\Di\times\N)$ is dense in $f^{-1}(U_i)$ by the corona theorem, each $\tau\circ g_i$ admits a continuous extension to a function on $f^{-1}(U_i)$ with range in $\tau(S_d)$, see \cite[Lm.\,8.1]{Br2}.
Applying to this extension the inverse map $\tau^{-1}:\tau(S_d)\rightarrow S_d$ we obtain that each $g_i$ admits a continuous extension $\hat g_i:f^{-1}(U_i)\rightarrow S_d$. Since the set $(f^{-1}(U_i)\cap (\Di\times\N))\cap (f^{-1}(U_j)\cap(\Di\times\N))$ is dense in the open set $f^{-1}(U_i)\cap f^{-1}(U_j)$ by the corona theorem, equation \eqref{e4.1} and continuity of the extensions imply that
\begin{equation}\label{e4.2}
\hat g_j(x)=\hat g_i(x)\cdot (f^*g_{ij})(x)\quad \text{for all}\quad x\in f^{-1}(U_i)\cap f^{-1}(U_j),\quad i,j\in I.
\end{equation}
The latter shows that the family $\{\hat g_i\}_{i\in I}$ determines a continuous section $\hat c$ of $f^*E$ over $\Pi$, as required.
\end{proof}
Using the lemma let us complete the proof of the proposition. 

We define
a map $f':\Pi\rightarrow \mathcal M'$ by the formula
\begin{equation}\label{e4.3}
f'(x)=\rho(\hat g_i(x))(s_i(f(x))),\quad x\in f^{-1}(U_i),\ i\in I.
\end{equation}
Since by \eqref{e4.0}, \eqref{e4.2}
\[
\begin{array}{r}
\rho(\hat g_i(x))(s_i(f(x)))=\rho(\hat g_i(x))\bigl(\rho(g_{ij}(f(x))(s_j(f(x)))\bigr)
=\rho\bigl(\hat g_i(x)(f^*g_{ij})(x)\bigr)(s_j(f(x))\\
\\
=\rho(\hat g_j(x))(s_j(f(x)))\quad {\rm for\ all}\quad x\in f^{-1}(U_i)\cap f^{-1}(U_j)\ne\emptyset,\quad\quad\
\end{array}
\]
the map $f'$ is well defined. Also, by the definition $r\circ f'=f$.

The proof of the proposition is complete.
\end{proof}
\noindent {\bf 5.2.} 
Let $K\subset \mathfrak M(H^\infty(\Di\times\N))$ be a holomorphically convex set 
and let $G$ be an abelian group. In this part we prove the following result.
\begin{Lm}\label{lem5.3}
The homomorphism
\begin{equation}\label{eq5.5}
H^1(\mathfrak M(H^\infty(\mathbb D\times\N)),G)\rightarrow H^1(\mathfrak M(A),G)
\end{equation}
induced by the restriction map to $K$ is surjective.
\end{Lm}
\begin{proof}
As was mentioned in the Introduction, Theorem \ref{te1} is equivalent to the fact that the above homomorphism is surjective for $G=\Z$. Hence, it is surjective for $G=\Z^k$, $k\in\N$, as well (because $H^1(X,\Z^k)=(H^1(X,\Z))^k$ for every topological space $X$).

In general, let $g=\{g_{ij}\in C(U_i\cap U_j,G)\}_{i,j\in I}$ be a $1$-cocycle defined on an open finite cover $\mathcal U=(U_i)_{i\in I}$ of $K$. Passing to a refinement of $\mathcal U$, if necessary, we may assume that each $g_{ij}$ is defined  on the closure  of $U_{i}\cap U_j$. Since the latter is a compact subset of $K$ and $G$ is considered with the discrete topology, the image of each $g_{ij}$ is finite. Then the subgroup of $G$ generated by elements of images of all $g_{ij}$ is a finitely generated abelian group $G'$. Hence, $g$ is a cocycle with values in $G'$. Next, group $G'$ being finitely generated is isomorphic to $G'_f\oplus G'_t$, where $G'_f\cong \Z^k$ for some $k\in \Z_+$ and $G'_f$ is the finite torsion subgroup of $G'$. Hence, we can write
$g=g_1+g_2$, where $g_1$ and $g_2$ are cocycles on $\mathcal U$ with values in $G_f'$ and $G_t'$, respectively. Cocyle $g_1$ determines an element $\{g_1\}\in H^1(K,G_f')\cong H^1(K,\Z^k)$. Thus as explained above, $\{g_1\}$ lies in the image of homomorphism \eqref{eq5.5} for $G=G_f'$. In turn, cocycle $g_2$ determines an element $\{g_2\}\in H^1(K,G_t')$. In particular, $\{g_2\}$ determines an isomorphism class of principal bundles on $K$ with finite fibre $G_t'$. According to Lemma \ref{lem4.2} every such bundle is trivial. This implies that  $\{g_2\}=0$; hence, we obtain that the cohomology class $\{g\}$ of the cocycle $g$ coincides with $\{g_1\}$. In particular, it lies in the image of homomorphism \eqref{eq5.5}, as required.
\end{proof}
\begin{proof}[{\bf 5.3.} Proof of Theorem \ref{teo1.3}]
Let $\varphi: H^\infty(\mathbb D\times\N)\rightarrow A$ be a dense image morphism of complex unital Banach algebras and let $\mathcal M\subset\Co^n$ be a connected complex submanifold and an Oka manifold. Assuming that a finite unbranched covering $\mathcal M'$ of $\mathcal M$ is $i$-simple for $i=1,2$, we have to prove that the image of  $(H^\infty(\mathbb D\times\N))^\mathcal M$ under $\varphi^n$ is a dense subset of $A_\mathcal M$. Retaining notation of Theorem \ref{te4.1} we assume that
$A^*\subset (H^\infty(\mathbb D\times\N))^*$ so that $\mathfrak M(A)$ is a compact holomorphically convex subset of $\mathfrak M(H^\infty(\mathbb D\times\N))$. Then since
 \[
 (H^\infty(\mathbb D\times\N))^\mathcal M=(H^\infty(\mathbb D\times\N))_\mathcal M
 \]
 because the algebra $H^\infty(\Di\times\N)$ is semisimple,
 due to this theorem, it suffices to prove that each $f\in C(\mathfrak M(A),\mathcal M)$ can be extended to a map $\tilde f\in C(\mathfrak M(H^\infty(\mathbb D\times\N)), \mathcal M)$.

Let $r: \mathcal M'\rightarrow \mathcal M$ be the covering map. According to Proposition \ref{prop5.1} there is a map $f'\in C(\mathfrak M(A),\mathcal M')$ such that $f=r\circ f'$.  Next, we use that ${\rm dim}\, \mathfrak M(H^\infty(\mathbb D\times\N))=2$, see \cite[Thm.\,2.6]{Br2}.
(Recall that for a normal space $X$, ${\rm dim}\,X \le n$ if every open cover of $X$ can be refined by an open cover whose order $\le n + 1$. If  ${\rm dim}\,X \le n$ and the statement  ${\rm dim}\,X \le n-1$ is false, then ${\rm dim}\,X = n$.) Since  $\mathcal M'$ is $1$-simple, the fundamental group $\pi_1(\mathcal M')$ is abelian. Then 
Due to Lemma \ref{lem5.3} the homomorphism
\[
 H^1(\mathfrak M(H^\infty(\mathbb D\times\N)),\pi_1(\mathcal M'))\rightarrow H^1(\mathfrak M(A),\pi_1(\mathcal M'))
\]
is surjective. This and \cite[(10.4)]{Hu2} with $n=1$ and $m=2$ imply that $f'$ can be extended to a map $\tilde f'\in C(\mathfrak M(H^\infty(\mathbb D\times\N)),\mathcal M')$.
We set $\tilde f:=r\circ\tilde f'\in C(\mathfrak M(H^\infty(\mathbb D\times\N)),\mathcal M)$. Then $\tilde f$ is the required extension of $f$.

The proof of the theorem is complete.
\end{proof}
\sect{Proofs of Theorems \ref{teo1.5} and \ref{teo1.9}}
\begin{proof}[Proof of Theorem \ref{teo1.5}]
Suppose $\Pi_c^k=\Pi_c^k[F]$, that is, it is determined by a family  $F=\{f_i\}_{1\le i\le k}\subset H^\infty$, see \eqref{e3.1}. The set $\mathcal D=\{F,\mathcal M, K, n, c, k,\delta,\varepsilon \}$ is called the {\em data}.  
Since $\Pi_c^k$ is the disjoint union of open sets biholomorphic to $\Di$, the corona theorem is valid for $H^\infty(\Pi_c^k)$.  Let $A\subset H^\infty(\Pi_{c;\delta}^k)$ be the uniform closure of the restriction of  $H^\infty$ to $\Pi_{c;\delta}^k$. According to
\cite[Cor.\,2.6]{S2} the restriction of $H^\infty(\Pi_c^k)$ to 
$\Pi_{c;\delta}^k$ lies in $A$. Hence, the transpose of the restriction embeds $\mathcal M(A)$ into $\mathfrak M(H^\infty(\Pi_c^k))$. Since by the hypotheses of the theorem $g\in (H^\infty(\Pi_c^k))^\mathcal M$, the latter implies that
$g|_{\Pi_{c;\delta}^k}\in A^\mathcal M\, (=A_\mathcal M)$.
Then according to Theorem \ref{teo1.3} given $\varepsilon>0$ there exists
$h\in (H^\infty)^\mathcal M$ (depending on $g$ and the data) such that
\begin{equation}\label{eq7.1}
 \|(h-g)|_{\Pi_{c;\delta}^k}\|_{(H^\infty(\Pi_{c;\delta}^k))^n}\le\varepsilon.
\end{equation}
By $\mathcal H_{g,\mathcal D}$ we denote the class of such maps $h$.
We have to prove that
\begin{equation}\label{e7.2}
C=C(\mathcal M, K,n, c,k,\delta,\varepsilon):=\sup_{F,g}\inf_{h\in \mathcal H_{g,\mathcal D}} 
\|h\|_{(H^\infty)^n}
\end{equation}
is finite. 

To this end, let  $\{\Pi_c^k[F_i]\}_{i\in\N}\subset \Di$, $F_i=\{f_{ji}\}_{j=1}^k\subset H^\infty$ and $\{g_i\}_{i\in\N}\subset H^\infty(\Pi_c^k,\mathcal M)$, $g_i(\Pi_c^k)\subset K$, be sequences  satisfying assumptions of the theorem such that
\begin{equation}
C=\lim_{i\rightarrow\infty}\inf_{h\in \mathcal H_{g_i,\mathcal D_i}}\|h\|_{(H^\infty)^n};
\end{equation}
here $\mathcal D_i:=\{F_i,\mathcal M, K, n,c, k, \delta,\varepsilon\}$.

As in the proof of Theorem \ref{lem3.1} we define $F=\{f_j\}_{1\le j\le k}\subset H^\infty(\Di\times\N)$, $\Pi_\nu[F]\subset \Di\times\N$ and $g\in \mathcal O(\Pi_1[F],\mathcal M)$, $g(\Pi_1[F])\subset K$, by the formulas
\begin{equation}\label{equ6.4}
f_j|_{\Di\times\{i\}}:=f_{ji},\ 1\le j\le k,\quad \Pi_\nu[F]\cap (\Di\times\{i\}):=\Pi_{c;\nu}^k[F_i],\quad g|_{\Pi_{c}^k[F_i]}:=g_i,\quad
i\in\N.
\end{equation}
Since $\Pi_1[F]$ is biholomorphic to $\Di\times\N$, the corona theorem is valid for $H^\infty(\Pi_1[F])$. This implies that $g\in (H^\infty(\Pi_1[F]))_\mathcal M$.

By $A_\delta [F]$ we denote the uniform closure of the restriction of $H^\infty(\Di\times\N)$ to $\Pi_\delta [F]$.
\begin{Lm}\label{lem6.1}
The restriction of $H^\infty(\Pi_1[F])$ to $\Pi_\delta [F]$ forms a dense subalgebra of $A_\delta [F]$.
\end{Lm}
\begin{proof}
Let $f\in H^\infty(\Pi_1[F])$. Then $f-c\in (H^\infty(\Pi_1[F]))^{-1}$, where $c:=\|f\|_{H^\infty(\Pi_1[F])}+1$, and 
\[
\|(f-c)^{\pm 1}\|_{H^\infty(\Pi_1[F])}\le 2c-1.
\]
Using this estimate and applying Theorem \ref{lem3.1} to each function
$(f-c)|_{\Pi_{c}^k[F_i])}$ we obtain that for every $\varepsilon>0$ there exist a constant $C$ and a function
$f_\varepsilon\in H^\infty(\Di\times\N)$ such that
\[
\|f_\varepsilon\|_{H^\infty(\Di\times\N)}\le C\quad {\rm and}\quad
\|(f_\varepsilon -(f-c))|_{\Pi_\delta [F]}\|_{H^\infty(\Pi_\delta [F])}\le\varepsilon .
\]
This shows that $f$ can be approximated uniformly on $\Pi_\delta [F]$ by restrictions of functions from $H^\infty(\Di\times\N)$, as required.
\end{proof}
From the lemma we obtain that the transpose of the restriction map $\Pi_1 [F]\rightarrow \Pi_\delta [F]$ embeds $\mathfrak M(A_\delta[F])$ into $\mathfrak M(H^\infty(\Pi_1 [F]))$. Then since
$g\in (H^\infty(\Pi_1 [F]))_{\mathcal M}$, the restriction
$g|_{\Pi_\delta [F]}\in (A_\delta [F])_{\mathcal M}$.
In particular, due to Theorem \ref{teo1.3} given $\varepsilon>0$ there exists
$f\in (H^\infty(\Di\times\N))^\mathcal M$
such that
\begin{equation}\label{eq6.5}
 \|(f-g)|_{\Pi_{\delta}[F]}\|_{(H^\infty(\Pi_{\delta}[F]))^n}\le\varepsilon.
\end{equation}
This implies that $C$ can be estimated from above by 
$\|f\|_{(H^\infty(\Di\times\N))^n}$ which completes the proof of the theorem.
\end{proof}
\begin{proof}[Proof of Theorem \ref{teo1.9}]
We retain notation of the proof of the previous theorem.

 Let $\Pi_c^k=\Pi_c^k[F]\subset\Di$. Our data in this case is the set $\mathcal D=\{F,\mathcal M, K,n,c,k,J,b,\delta\}$.
Let $r:H^\infty \rightarrow C({\rm hull}\, J)$,   $r(f):=\hat f|_{{\rm hull}\, J}$, be the restriction homomorphism and let 
$I={\rm ker}\, r$. Then $I\subset H^\infty$ is a proper closed ideal containing $J$. We naturally identify the complex unital Banach algebra $A_I:=H^\infty/I$ with the algebra $r(H^\infty)$ equipped with the quotient norm.
Since ${\rm hull}\, I={\rm hull}\, J$, the transpose of $r$ maps the maximal ideal space $\mathfrak M(A_I)$ homeomorphically onto ${\rm hull}\, J$. We identify these two spaces so that each homomorphism from $\mathfrak M(A_I)$ is the evaluation homomorphism at a point of ${\rm hull}\, J$.  By the hypothesis of the theorem $r^n(g)\in (A_I)_\mathcal M$, see \eqref{equ1.11}; hence, Theorem \ref{teo1.3} implies that there exists $h\in (H^\infty)^\mathcal M$ (depending on $g$ and the data) such that 
\begin{equation}\label{equ6.6}
r^n(h)=r^n(g).
\end{equation}
By $\mathcal H_{g,\mathcal D}$ we denote the class of such maps $h$.
We have to prove that
\begin{equation}\label{equ6.7}
C=C(\mathcal M, K,n, b, c,k,\delta):=\sup_{F,J,g}\inf_{h\in \mathcal H_{g,\mathcal D}} 
\|h\|_{(H^\infty)^n}
\end{equation}
is finite. 

As before, let  $\{\Pi_c^k[F_i]\}_{i\in\N}\subset \Di$, $J_i\subset H^\infty$,  $\{g_i\}_{i\in\N}\subset (H^\infty)^n$,  be sequences  satisfying assumptions of the theorem such that
\begin{equation}\label{equ6.8}
C=\lim_{i\rightarrow\infty}\inf_{h\in \mathcal H_{g_i,\mathcal D_i}}\|h\|_{(H^\infty)^n};
\end{equation}
here $\mathcal D_i:=\{F_i,\mathcal M, K, n,c, k, J_i, b,\delta\}$.

By $r_i:H^\infty\rightarrow C({\rm hull}\, J_i)$ we denote the corresponding restriction homomorphisms, by $I_i$ their kernels, and by $A_{I_i}$ the corresponding quotient algebras. The hypotheses of the theorem provide also some maps $f_i\in\mathcal O(\Pi_c^k[F_i],\mathcal M)$ with images in $K$ such that
\begin{equation}\label{eq6.9}
\hat f_i|_{{\rm hull}\, I_i}=\hat g_i|_{{\rm hull}\, I_i}.
\end{equation}
As in the proof of the previous theorem, see \eqref{equ6.4}, we define $F\subset H^\infty(\Di\times\N)$, $\Pi_\nu[F]\subset \Di\times\N$, $g\in H^\infty(\Di\times N)$, and $f\in \mathcal O(\Pi_1[F],\mathcal M)$, $f(\Pi_1[F])\subset K$, by the formulas
\begin{equation}
F|_{\Di\times\{i\}}:=F_i,\ \, \Pi_\nu[F]\cap (\Di\times\{i\}):=\Pi_{c;\nu}^k[F_i],\ \,  g|_{\Di\times\{i\}}:=g_i,\ \, f|_{\Pi_{c}^k[F_i]}:=f_i,\ \,
i\in\N.
\end{equation}
Since $\Pi_1[F]$ is biholomorphic to $\Di\times\N$, the corona theorem is valid for $H^\infty(\Pi_1[F])$. This implies that $f\in (H^\infty(\Pi_1[F]))_\mathcal M$. Moreover, according to \cite[Lm.\,8.1]{Br2}, $f$ admits a continuous extension $\hat f$ to an open polyhedron $\hat \Pi_1[\hat F]\subset\mathfrak M(H^\infty(\Di\times\N))$ determined by extension (denoted by $\hat F$) of the family $F$ to $\mathfrak M(H^\infty(\Di\times\N))$ by the Gelfand transform (cf. \eqref{e3.10}). Also, $\Pi_1(F)$ is an open dense subset of $\hat \Pi_1[\hat F]$ by the corona theorem for $H^\infty(\Di\times\N)$. This implies that $\hat f$ maps $\hat \Pi_1[\hat F]$ into $K$.

\medskip
Next, by the definition, 
\[
\Di\times\N\subset\mathfrak M(H^\infty)\times\N\subset\mathfrak M(H^\infty(\Di\times\N)).
\]
Let 
\begin{equation}\label{equ6.11}
r:H^\infty(\Di\times\N)\rightarrow C_b\left(\bigsqcup_{i\in\N}{\rm hull}\, J_i\right),\quad r(f)|_{{\rm hull}\, J_i}:=r_i(f|_{\Di\times\{i\}}),\quad i\in\N.
\end{equation}
Clearly, $r$ is a well-defined morphism of complex Banach algebras of norm $\le 1$ (Here $C_b(X)\subset C(X)$ is the Banach algebra of bounded continuous functions on $X$ equipped with supremum norm.) We set $I:={\rm ker}\, r$. Then
\begin{equation}\label{equ6.12}
u\in I\Longleftrightarrow u|_{\Di\times\{i\}}\in I_i\quad \forall\,  i\in\N.
\end{equation}
We define $A_I:=H^\infty(\Di\times\N)/I$ and naturally identify this algebra with $r(H^\infty(\Di\times\N))$ equipped with the quotient norm. Now the transpose of $r$ maps the maximal ideal space $\mathfrak M(A_I)$ homeomorphically onto the set
\[ 
{\rm hull}\, I:=\{\xi\in \mathfrak M(H^\infty(\Di\times\N))\, :\, \hat u(\xi)=0\quad \forall\, u\in I\}.
\]
\begin{Lm}\label{lem6.2}
It is true that
\[
{\rm hull}\, I\subset \hat \Pi_1[\hat F].
\]
\end{Lm}
\begin{proof}
Let $\bar A_I$ be the uniform closure of $A_I$ in $C_b(X)$, $X:=\underset{i\in\N}\sqcup\, {\rm hull}\, J_i$. Then the transpose of the embedding $A_I\hookrightarrow \bar A_I$ maps $\mathfrak M(\bar A_I)$ homeomorphically onto $\mathfrak M(A_I)={\rm hull}\, I$.

Let $\hat F=\{\hat f_j\}_{1\le j\le k}$, $f_j\in H^\infty(D\times\N)$. Suppose on the contrary that there is $\xi\in  {\rm hull}\, I\setminus \Pi_1[\hat F]$. 
Then since each $\xi\in\mathfrak M(\bar A_I)$ has norm $\le 1$, by assumption \eqref{subset} of the theorem we obtain
\[
1\le \max_{1\le j\le k}|\hat f_j(\xi)|\le \max_{1\le j\le k}\sup_{X}|\hat f_j |\le\delta<1.
\]
This contradiction proves the required implication.
\end{proof}
Lemma \ref{lem6.2} and the fact that $\hat f$ maps $\hat \Pi_1[\hat F]$ into $K$ and  $\hat f|_X=r^n(g)$, see \eqref{eq6.9},
show that $r^n(g)\in (A_I)_\mathcal M\, (=(A_I)^\mathcal M)$.

Finally, applying Theorem \ref{teo1.3} to the epimorphism
$r:H^\infty(\Di\times\N)\rightarrow A_I$ and the element $r^n(g)$ we find an element $h\in H^\infty(\Di\times\N)^\mathcal M$ such that
\[
r^n(h)=r^n(g).
\]
Then by our definition, see \eqref{equ6.8},
$h|_{\Di\times\{i\}}\in\mathcal H_{g_i,\mathcal D_i}$ and  
\[
C\le \sup_{i\in\N}\|h|_{\Di\times\{i\}}\|_{(H^\infty)^n}=:\|h\|_{(H^\infty(\Di\times\N))^n}<\infty.
\]
This completes the proof of the theorem.
\end{proof}

\end{document}